\documentclass{article}
\usepackage{latexsym} 
\usepackage[english]{babel} 
\usepackage{epsfig} 
\usepackage{amsmath} 
\usepackage{amsfonts} 
\usepackage{amssymb}
\usepackage{graphicx} 
\usepackage{tikz}
\usetikzlibrary{arrows}
\usetikzlibrary{calc}
\usetikzlibrary{decorations.pathmorphing}

\newcommand{\sgn}{\operatorname{sgn}}
\newcommand{\ddiv}{\operatorname{div}}
\newcommand{\cpg}{\mathop{\Box}} 

\tikzset{snake it/.style={decorate, decoration={snake, amplitude=0.9, segment length=3}}}

\newlength{\pinner}
\newlength{\pouter}
\newlength{\thickedge}
\newlength{\thinedge}
\definecolor{farcolor}{rgb}{0.8,0.2,0.6}
\definecolor{nearcolor}{rgb}{0,0.7,0}
\definecolor{zigzagcolor}{rgb}{0.3,0.8,1}
\definecolor{unbalancedcolor}{gray}{0.8}
\tikzstyle{far}=[black, line width = \thickedge, dashed, shorten >=2pt, shorten <=2pt]
\tikzstyle{near}=[black, line width = \thickedge, shorten >=2pt, shorten <=2pt]
\tikzstyle{zigzag}=[black!80, line width = \thinedge, shorten >=4pt, shorten <=4pt, snake it]

\tikzstyle{vertex}=[circle, draw=black!70, minimum width=2.5, inner sep=0pt,fill=black!20]
\tikzstyle{blackvertex}=[circle, draw=black, minimum width=4, inner sep=0pt, fill=black]
\tikzstyle{defaultedge}=[line width = 1pt, color = black!30]
\tikzstyle{defaultdiredge}=[line width = 1pt, shorten >=1pt, color = black!30]
\tikzstyle{dte}=[line width = 1pt, color = black!10,shorten >=2pt, shorten <=2pt]
\tikzstyle{unbalanced}=[fill=unbalancedcolor]

\newcommand{\petersencoords}{
\coordinate (p0) at (0:\pouter);
\coordinate (p1) at (72:\pouter);
\coordinate (p2) at (144:\pouter);
\coordinate (p3) at (216:\pouter);
\coordinate (p4) at (288:\pouter);
\coordinate (p5) at (0:\pinner);
\coordinate (p6) at (72:\pinner);
\coordinate (p7) at (144:\pinner);
\coordinate (p8) at (216:\pinner);
\coordinate (p9) at (288:\pinner);
}
\newcommand{\petersennodes}{
\foreach \x in {p0, p1, p2, p3, p4, p5, p6, p7, p8, p9}{
  \node[vertex] at (\x) {};
}
}

\newcommand{\redpetersen}{
\draw[defaultedge] (-1,0) --node[black, above] {${e}$} (1,0); 
\draw[defaultdiredge,->] (-1,0) -> (-1,1.4) node[black, pos=0.9, left]{$v_3$};
\draw[defaultdiredge,->] (1,1.4) -> (1,0) node[black, pos=0.1, right]{$v_5$};
\draw[defaultdiredge,->] (0,2) -> (1,1.4);
\draw[defaultdiredge,->] (-1,1.4) -> (0,2) node[black, pos=0.95, above left]{$v_4$};
\draw[defaultdiredge,->] (-1,0) -> (-2.4,0) node[black, pos=0.9, above]{$v_2$} node[black, pos=0, above left]{$v_1$};
\draw[defaultdiredge,->] (1,0) -> (2.4,0) node[black, pos=0.9, above]{$v_7$} node[black, pos=0, above right]{$v_6$};

\node[vertex] at (-1,0) {};
\node[vertex] at (1,0) {};
\node[vertex] at (-1,1.4) {};
\node[vertex] at (1,1.4) {};
\node[vertex] at (0,2) {};
\node[vertex] at (-2.4,0) {};
\node[vertex] at (2.4,0) {};
}

\newcommand{\petersenedges}{
  \draw[defaultedge] (p0) -- (p1);
  \draw[defaultedge] (p1) -- (p2);
  \draw[defaultedge] (p2) -- (p3);
  \draw[defaultedge] (p3) -- (p4);
  \draw[defaultedge] (p4) -- (p0);
  \draw[defaultedge] (p0) -- (p5);
  \draw[defaultedge] (p1) -- (p6);
  \draw[defaultedge] (p2) -- (p7);
  \draw[defaultedge] (p3) -- (p8);
  \draw[defaultedge] (p4) -- (p9);
  \draw[defaultedge] (p5) -- (p7);
  \draw[defaultedge] (p7) -- (p9);
  \draw[defaultedge] (p9) -- (p6);
  \draw[defaultedge] (p6) -- (p8);
  \draw[defaultedge] (p8) -- (p5);
}

\newcommand{\petersenlabelededges}{
  \draw[defaultedge] (p0) --node[black, above right] {${e}_3$} (p1);
  \draw[defaultedge] (p1) --node[black, above left] {${e}_4$} (p2);
  \draw[defaultedge] (p2) --node[black, left] {${e}_0$} (p3);
  \draw[defaultedge] (p3) --node[black, below] {${e}_1$} (p4);
  \draw[defaultedge] (p4) --node[black, below right] {${e}_2$} (p0);
  \draw[defaultedge] (p0) --node[black, below] {${f}_0$} (p5);
  \draw[defaultedge] (p1) --node[black, left] {${f}_1$} (p6);
  \draw[defaultedge] (p2) --node[black, below] {${f}_2$} (p7);
  \draw[defaultedge] (p3) --node[black, right] {${f}_3$} (p8);
  \draw[defaultedge] (p4) --node[black, right] {${f}_4$} (p9);
  \draw[defaultedge] (p5) -- (p7) node[black, pos=0.8, above]{$e'_1$};
  \draw[defaultedge] (p7) -- (p9) node[black, pos=0.88, left]{$e'_3$};
  \draw[defaultedge] (p9) -- (p6) node[black, pos=0.85, right]{$e'_0$};
  \draw[defaultedge] (p6) -- (p8) node[black, pos=0.8, left]{$e'_2$};
  \draw[defaultedge] (p8) -- (p5) node[black, pos=0.85, below]{$e'_4$};
}

\newcommand{\fan}{
\foreach \x in {1,...,8}{
  \draw[defaultedge] (\x,0) -- +(1,0);
}
\foreach \x in {1,...,9}{
  \draw[defaultedge] (\x,0) -- (5,2);
}
\foreach \x in {1,...,9}{
  \node[vertex] at (\x,0) {}; 
}  
\node[vertex] at (5,2) {}; 
\foreach \x in {1,...,8}{
  \draw[zigzag] (\x,0) -- +(1,0); 
}
\node at (5.5,2.1) {$v_0$};
\foreach \y in {1,...,9}{
  \node at (\y,-0.3) {$v_\y$};
}
}

\newcommand{\ladder}[1]{
\foreach \x in {1,...,#1}{
  \draw[defaultedge] (\x-1,0) -- +(1,0);
  \draw[defaultedge] (\x-1,1) -- +(1,0);
  }
\foreach \x in {0,...,#1}{
  \draw[defaultedge] (\x,0) -- +(0,1);
  }
\foreach \x in {0,...,#1}{
  \node[vertex] at (\x,0) {};
  \node[vertex] at (\x,1) {};
  }   
}

\newcommand{\labeledladder}[1]{   
\foreach \x in {0,...,#1}{
  \draw[defaultedge] (\x-1,1) --node[black, above] {$\overline{e}_\x$} +(1,0);
  \draw[defaultedge] (\x-1,0) --node[black, below] {$\underline{e}_\x$} +(1,0);
  }
\foreach \x in {0,...,#1}{
  \draw[defaultedge] (\x-1,0) -- node[black, right] {$e'_\x$} +(0,1);
  }
\foreach \x in {#1,...,#1}{  
  \draw[defaultedge] (\x,0) -- node[black, right] {$e'_4$} +(0,1); 
}  
\foreach \x in {0,...,#1}{
  \node at (\x-1,-0.2) {$u_\x$};
  \node at (\x-1,1.2) {$v_\x$};
  } 
\foreach \x in {#1,...,#1}{  
  \node at (\x,-0.2) {$u_4$};
  \node at (\x,1.2) {$v_4$};
} 
\foreach \x in {0,...,#1}{
  \node at (\x-0.2,0.2) {$\sigma_\x$};
  }
\foreach \x in {0,...,#1}{
  \node[vertex] at (\x,0) {};
  \node[vertex] at (\x,1) {};
  } 
\node[vertex] at (-1,0) {};
\node[vertex] at (-1,1) {};    
}

\newcommand{\labeledprism}[1]{   
\foreach \x in {0,...,#1}{
  \draw[defaultedge] (\x-1,1) --node[black, above] {$\overline{e}_\x$} +(1,0);
  \draw[defaultedge] (\x-1,0) --node[black, below] {$\underline{e}_\x$} +(1,0);
  }
\foreach \x in {0,...,#1}{
  \draw[defaultedge] (\x-1,0) -- node[black, right] {$e'_\x$} +(0,1);
  }
\foreach \x in {#1,...,#1}{  
  \draw[defaultedge] (\x,0) -- node[black, right] {$e'_0$} +(0,1);  
}  
\foreach \x in {0,...,#1}{
  \node at (\x-1,-0.2) {$u_\x$};
  \node at (\x-1,1.2) {$v_\x$};
  } 
\foreach \x in {#1,...,#1}{  
  \node at (\x,-0.2) {$u_0$};
  \node at (\x,1.2) {$v_0$};
} 
\foreach \x in {0,...,#1}{
  \node at (\x-0.2,0.2) {$\sigma_\x$};
  }
\foreach \x in {0,...,#1}{
  \node[vertex] at (\x,0) {};
  \node[vertex] at (\x,1) {};
  } 
\node[vertex] at (-1,0) {};
\node[vertex] at (-1,1) {};    
}

\newcommand{\unbalancedsquare}[1]{
  \fill[unbalanced] (#1,0) -- ++(-1,0) -- ++(0,1) -- ++ (1,0) -- cycle 
}

\newcommand{\ladderplus}[1]{
\foreach \x in {-1,...,#1}{
  \draw[defaultedge] (\x,0) -- +(1,0);
  \draw[defaultedge] (\x,1) -- +(1,0);
  }
\foreach \x in {0,...,#1}{
  \draw[defaultedge] (\x,0) -- +(0,1);
  }
\foreach \x in {0,...,#1}{
  \node[vertex] at (\x,0) {};
  \node[vertex] at (\x,1) {};
  }   
}

\newtheorem{thm}{Theorem}[section]

\newtheorem{lem}[thm]{Lemma}
\newtheorem{cor}[thm]{Corollary}

\begin{document} 

\title{{\bf Threshold Colorings of Prisms \\ and the Petersen Graph}}

\author{Ga\v{s}per Fijav\v{z} \and Matthias Kriesell}

\maketitle



\begin{abstract}
  Let $G$ be a graph, $r \geq t$ integers, and $N \subseteq E(G)$.
  An {\em $(r,t)$-threshold-coloring of $G$ with respect to $N$} is a mapping $c: V(G) \rightarrow \{0,\ldots,r-1\}$ such that
  $|c(u)-c(v)| \leq t$ for every $uv \in N$ and $|c(u)-c(v)|>t$ for every $uv \in E(G) \setminus N$.
  A graph is {\em total threshold colorable} if there exist integers $r,t$ such that for every $N \subseteq E(G)$,
  $G$ admits an $(r,t)$-threshold-coloring with respect to $N$.
  We show that every prism is total threshold colorable, and that the {\sc Petersen} graph is total threshold colorable.
  In contrast to this fact we show that {\sc M\"obius} ladders are {\em not} total threshold colorable, from which
  it follows that there is no characterization of being total threshold colorable in terms of a finite set of forbidden subgraphs.
\end{abstract}



\section{Introduction}

All graphs considered here are supposed to be finite, simple, and undirected.
A {\em near-far-labeling} of a graph $G=(V,E)$ is a pair $(N,F)$ with $N \subseteq E(G)$ and $F=E(G) \setminus N$,
where the edges from $N$ will be called {\em near edges} and those of $F$ are the {\em far edges}.
For integers $r \geq t$, an {\em $(r,t)$-threshold-coloring with respect to the near-far-labeling $(N,F)$} is a mapping
$c : V(G) \rightarrow \{0,\ldots,r-1\}$ such that for every near edge $uv \in N$ we have $|c(u) - c(v)| \leq t$,
and for every far edge  $uv \in F$ we have $|c(u) - c(v)| > t$.
In other words, end vertices of a near edge receive {\em near} colors in terms of the absolute value of their difference, end vertices of far edges receive colors that are far apart.
The notion of near and far is measured with integer $t$, which we call the {\em threshold} of the coloring,
and colors stem from the set of the first $r$ non-negative integers;  $r$ is the {\em range} of the coloring. If more convenient,
we shall use alternative sets of $r$ consecutive integers.
A graph $G$ is {\em $(r,t)$-total-threshold-colorable} if there exists an $(r,t)$-threshold-coloring of $G$ with respect to {\em every} near-far-labeling.
Disregarding parameters   we say that $G$ is {\em total threshold colorable} if it is $(r,t)$-total-threshold-colorable for some integers $r \geq t$.

These concepts have been introduced in \cite{AlamChaplickFijavzKaufmannKobourovPupyrev2013},
where various complexity issues and relations to other coloring problems have been discussed
(like, for example, to questions on distance graphs and threshold graphs as introduced in \cite{HammerPeledSun1990} and \cite{MahadevPeled1995}, respectively).

As an application, one can infer that if some planar graph has a unit-cube contact representation and is, at the same time, total threshold colorable,
then all its subgraphs are unit-cube contact representable, too \cite{AlamChaplickFijavzKaufmannKobourovPupyrev2013}.
Another (straightforward) application is a puzzle called {\em Happy Edges}, played on a square grid with solid and dotted edges where one has to put
integer values on the vertices such that they differ by at most one along a solid edge, and by more than one along a dotted edge
\cite{AlamKobourovPupyrevToeniskoetter2014}.

One of the interesting open problems is whether finite square grids are total threshold colorable;
in \cite{AlamKobourovPupyrevToeniskoetter2014} it has been proved that for every $r>0$ there exists a square grid (of large size) which is not $(r,t)$-total-threshold-colorable for all $t$,
so that one necessarily has to increase the range when looking for threshold colorings of larger and larger grids.
Several classes of planar graphs such as cycles, trees, or, more general, planar graphs of girth at least $10$
turned out to be $(r,t)$-threshold-colorable for small values of $r,t$ (independent of their size).
However, all these graphs have vertices of degree at most $2$ (and are planar).
Here we study two infinite classes of $3$-connected graphs (indeed, classes of $3$-connected near-far-labeled graphs),
namely prisms, which turn out to be total threshold colorable, and {\sc M\"obius} ladders, which are not ---
a perhaps surprising fact, as, locally, a {\sc M\"obius} ladder is not distinguishable from a prism.
Another interesting problem is whether there is a good characterization of (total-) threshold-colorable graphs,
and one can infer from the result on {\sc M\"obius} ladders that, at least, one cannot characterize them in terms of finitely many forbidden subgraphs.

For two graphs $G,H$ let 
$G \cpg H$ 
denote their {\em cartesian product},
whose vertices are the pairs from $V(G) \times V(H)$, and where two distinct vertices $(v,w),(x,y)$ are adjacent if and only if either $v=x$ and $wy \in E(H)$,
or $w=y$ and $vx \in E(G)$. The subgraphs induced in $G \cpg H$ by sets of the form $V(G) \times \{y\}$ with $y \in V(H)$ or $\{x\} \times V(H)$ with $x \in V(G)$
are called the {\em $G$-fibres} or {\em $H$-fibres}, respectively. 
The graph $L_n:=P_n \cpg K_2$, where $P_n$ is the path of length $n$, is called a {\em ladder}. For $n \geq 3$, 
$L_n$ contains exactly two edges
whose end vertices both have degree $2$, and
by identifying these edges (and their endpoints) we get either $C_n \cpg K_2$, that is,
the {\em prism over the cycle $C_n$ of length $n$}, or a non-planar graph $M_n$, usually called the {\em M\"obius ladder of length $n$}.

Whereas it follows easily from a generalization of the methods in \cite{AlamChaplickFijavzKaufmannKobourovPupyrev2013}
that the ladders are $(5,1)$-threshold-colorable with respect to any near-far-labeling,
it took us much more effort to prove that prisms are total threshold colorable (Section \ref{SPrisms}).

\begin{thm} 
  \label{prismain}
  For every $n \geq 3$, the prism $C_n \cpg K_2$ is $(31,4)$-total-threshold-colorable.
\end{thm}

For {\sc M\"obius} ladders, the situation is different.

\begin{thm}
  \label{moebiusmain}
  For every $n \geq 3$, the {\sc M\"obius} ladder $M_n$ is not total threshold colorable.
\end{thm}

As a corollary, we get (almost immediately) the following:

\begin{cor}
  \label{cormoebius}
  There is no finite set $\mathfrak F$ of graphs such that
  a graph is total threshold colorable if and only if it contains no [induced] subgraph isomorphic to one of those in $\mathfrak F$.
  In other words, the property of being total threshold colorable is not characterizable in terms of a finite set of forbidden [induced] subgraphs.
\end{cor}

Finally, we look at the {\sc Petersen} graph, which is often seen as a problematic candidate in terms of colorings, in particular when coloring its edges.
In the threshold variant of vertex coloring we will ultimately succeed (without putting too much effort on optimal parameters).
Furthermore, the {\sc Petersen} graph is an example of a $3$-connected nonplanar cubic total threshold colorable graph (and the only one we know so far).

\begin{thm}
  \label{petersenmain}
  The {\sc Petersen} graph is total threshold colorable.
\end{thm}

The paper is organized accordingly: In Section \ref{STools} we introduce some basic tools,
in Section \ref{SLadders} we discuss ladders and {\sc M\"obius} ladders, 
in Section \ref{SPrisms} we give a proof of Theorem \ref{prismain}, and in the final section we treat the {\sc Petersen} graph.

\section{Toolbox}
\label{STools}

This section contains basic tools we need to prove our main results.  Let us start with the following lemma.

\begin{lem}
\label{init}
  Suppose that $c$ is an $(r,t)$-threshold-coloring of a graph $G$ with respect to the near-far-labeling $(N,F)$. 
      If $wxyzw$ is a cycle of length $4$ and $wx,yz$ are far edges and $xy,zw$ are near edges
      then $c(w)>c(x)$ implies 
$c(z)>c(y)$.  
\end{lem}

{\bf Proof.}
We infer $c(w)>c(x)+t$ from the assumption that $c(w)>c(x)$ and $wx$ is a far edge.
Since $yz$ is a far edge, either $c(z)>c(y)+t$ or $c(y)>c(z)+t$. In the first case we are done, so assume, to the contrary, that 
$c(y)>c(z)+t$. 
Consequently, $c(y)-c(x)=(c(y)-c(z))+(c(z)-c(w))+(c(w)-c(x))>t-t+t$, contradiction. 
\hspace*{\fill}$\Box$
 
 From Lemma \ref{init} it follows easily, that the complete graph $G:=K_4$ on four vertices is not total threshold colorable:
 Let $wxyzw$ be a spanning $4$-cycle, let $F=\{wx,yz\}$, and let $N=E(G) \setminus F$. Without loss of generality, $c(w)>c(x)$ so that $c(y)>c(z)$ by Lemma \ref{init};
 however, we may apply Lemma \ref{init} to the $4$-cycle $wxzyw$, implying $c(z)<c(y)$, contradiction.  A generalization of this argument will show that
 {\sc M\"obius} ladders are not total threshold colorable.

A {\em fan} $F_n$ is a graph obtained by adding a universal (or apex) vertex $v_0$ to the $n$-vertex path $v_1 v_2 \ldots v_n$.
We shall call this $n$-path the {\em spine} of $F_n$.
It was shown in \cite{AlamChaplickFijavzKaufmannKobourovPupyrev2013} that $F_n$ is $(5,1)$-threshold-colorable with respect to an arbitary near-far-labeling $(N,F)$.
We shall repeat the argument here both for clarity, and also for generalization purposes.

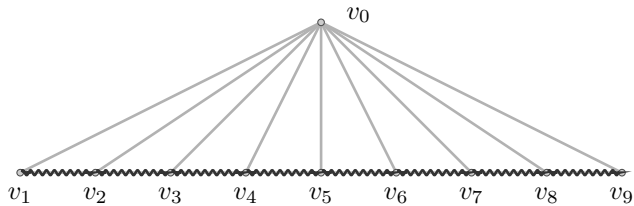
\begin{figure}[h]
\begin{center}
\begin{tikzpicture}[-]
\fan
\end{tikzpicture}
\end{center}
\caption{The fan $F_9$, with its spine shown in zigzag.}
\label{fig:fan}
\end{figure}

\begin{sloppy}
We shall use colors $0,\pm 1, \pm 2$ and inductively construct a $(5,1)$-threshold-co\-lo\-ring $c$.
Let us set $c(v_0)=0$.
For the remaining vertices we shall use colors $\pm 1,\pm 2$. 
Let us set $c(v_1)=1$ or $c(v_1)=2$ depending on whether $v_0 v_1$ is a near or a far edge, respectively.
Assuming that the colors $c(v_1),\ldots,c(v_i)$ are set, let us choose the color $c(v_{i+1})$ of the next vertex along the spine according to the following two criteria: 
\end{sloppy}

\begin{itemize}
  \item[(i)] the absolute value $|c(v_{i+1})|$ is either $1$ or $2$, depending on whether $v_0 v_{i+1}$ is near or far, respectively, and 
  \item[(ii)] the sign of $c(v_{i+1})$ is the same as $\sgn(c(v_i))$ if and only if $v_i v_{i+1}$ is a near edge.
\end{itemize} 

It is easy to verify that $c$ is a $(5,1)$-coloring.
In \cite{AlamKobourovPupyrevToeniskoetter2014}, the above argument has been extended to graphs which can be decomposed into a $2$-independent set
(that is, a set of vertices which are pairwise at distance at least $2$) and a collection of disjoint trees.
Here we shall extend the fan coloring argument a couple of steps further, as follows.

\begin{lem}
  \label{zigzag}
  Let $G=(V,E)$ be a graph and $(N,F)$ be a near-far-labeling. Suppose that there exists a bipartition of $V(G)$ into sets $A,B$ and a set $M \subseteq E(G)$ with the following properties.
\begin{itemize}
  \item[(i)]
    $E(G[A]) \subseteq N$,
  \item[(ii)]
    every cycle of $G[B]$ contains an even number of edges from $F$,
  \item[(iii)]
    $M \subseteq E(G)$ is an induced matching (that is $V(e),V(f)$ are at distance at least $2$ for $e \not= f$ from $M$)
    such that every edge from $M$ connects a vertex from $A$ to one from $B$, and
  \item[(iv)]
    for every pair of edges $ab,a'b \in E(G) \setminus M$ with $a,a' \in A$ and $b \in B$, we have either $ab,a'b \in N$ or $ab,a'b \in F$
\end{itemize}
Then there exists a $(13,4)$-threshold-coloring of $G$ with respect to $(N,F)$.
Moreover, if $M = \emptyset$ then there exists a $(5,1)$-threshold-coloring of $G$.
\end{lem}

{\bf Proof.}
Suppose that $G,(N,F),A,B,M$ satisfy all four conditions of the lemma.
The edges of $M$ can be interpreted as obstructions to condition (iv).
We thus may assume that all edges from $M$ do violate condition (iv), that is, for every edge $ab \in M$ with $a \in A,b \in B$ there exists an edge $a'b \in E(G)$ such that
$a' \in A$ and exactly one of $ab,a'b$ is from $N$. In fact, if $a''b$ is another such edge then both $a'b,a''b$ are from the same set $N$ or $F$ by (iv) as they are not in $M$.
Let us first relabel the edges from $M$, that is, we define an alternative near-far-labeling by
\[ (N',F'):=(N \triangle M,F \triangle M), \]
where $\triangle$ stands for the symmetric difference of sets.
It is easy to check that $G,(N',F'),A,B,M':=\emptyset$ satisfy all four conditions of the Lemma, too.
Observe that $(N',F')=(N,F)$ if (and only if) $M=\emptyset$.

We shall first construct a $(5,1)$-threshold-coloring $c'$ of $G$ with respect to $(N',$ $F')$. 
We will perturb $c'$ later in order to obtain a $(13,4)$-threshold-coloring of $G$ with respect to $(N,F)$. 

Let us first set $c'(v) = 0$ for every $v \in A$. This is a valid assignment as every edge between vertices from $A$ is a near edge. 

Let $T$ be a maximal spanning forest of $G[B]$ and let $R = E(G[B]) \setminus E(T)$.
Let us choose a root in every component of $T$, and orient the edges of $T$ away from respective roots.
We shall color the vertices of $T$ inductively along edge orientations and finally argue that the conditions imposed by labels of edges in $R$ are satisfied, too.

If $v$ is a root, then set $c'(v) = 1$ if all the edges connecting $v$ to vertices from $A$ are in $N'$, and set $c'(v)=2$ if they are all in $F'$ (other cases do not occur by (iv)).
Assume that $c'$ is partially determined and choose an edge $bd \in E(T)$ so that $c'(b)$ is determined and $d$ is yet uncolored.
We define $c'(d)$ by choosing its sign and absolute value.
If $bd \in N'$ then we choose $\sgn(c'(d)) = \sgn(c'(b))$, otherwise we choose $\sgn(c'(d)) = -\sgn(c'(b))$. 
If $d$ is incident with a near edge $ad$ where $a \in A$, then we choose $|c'(d)|=1$, otherwise we choose $|c'(d)|=2$.
Note that if $d$ is adjacent to a vertex from $A$, then by (iv) the edges linking $d$ to $A$ all have {\em the same} label.

Observe that for every edge $bd \in E(T)$ the colors $c'(b)$ and $c'(d)$ have he same sign if and only if $bd \in N$.
Let $b d \in R$. By (ii) the fundamental cycle in $T + bd$ contains an even number of edges from $F$.
This implies that the number of sign changes along the $b,d$-path in $T$ is even if and only if $bd$ is from $N'$.
Now, every pair of nonzero colors of the same sign are at most $1$ apart, and every pair of colors of opposite signs are more than $1$ apart.
Consequently, $c'$ is a $(5,1)$-threshold-coloring of $G$ with respect to $(N',F')$, and it proves the moreover-part of the Lemma.

Now let $3c'$ denote the coloring with $-6,-3,0,3,6$ obtained by multiplying the colors of $c'$ by $3$.
Note that $3c'$ is a $(13,3)$- and $(13,4)$- and also a $(13,5)$-threshold-coloring of $G$ with respect to $(N',F')$.
We go for a $(13,4)$-threshold-coloring,
which we have ready in hands  if $M=\emptyset$. Hence we may assume that $M \not= \emptyset$.
Then $3c'$ is not a $(13,4)$-threshold-coloring of $G$ with respect to $(N,F)$
because the span condition is violated on every single edge of $M$. We need to recolor endvertices of $F$.

Consider $ab \in M$ with $a \in A$, $b \in B$. Then $c'(a)=0$, and we may assume that $c'(b)>0$ first.
If $ab \in N$ then $ab \in F'$, so that $c'(b)=2$, and we set $c(b)=5 = 3 c'(b)-1$ and $c(a) = +1 = 3c'(a)+1$. 
If $ab \in F$ then $ab \in N'$, so that $c'(b)=1$, and we set $c(b)= 4 = 3 c'(b)+1$ and $c(a) = -1 = 3c'(a)-1$.
If $c'(b)<0$ then we do the same assignments with opposite signs.
For every vertex $x$ not incident with any edge from $M$ we set $c(x)=3c'(x)$. 
For every $e \not\in M$ we have thus altered the color of at most one endvertex (compared to $3c'$),
as $M$ is an induced matching, by $\pm1$ by construction.
Consequently, $c$ a $(13,4)$-threshold-coloring of $G$ with respect to $(N,F)$.  
\hspace*{\fill}$\Box$

The conditions of Lemma~\ref{zigzag} are technical and appear difficult to verify.
However, for most applications, we will be in a rather simple scenario.
As a rule of thumb almost every vertex of $B$ will be adjacent to at most one vertex of $A$. This will ensure (iv) everywhere except at a very small subset of $B$.
Similarly, the graph $G[A]$ will typically induce at most two (near) edges, and the graph $G[B]$ will either be a single path, a union of two paths, or a
unicyclic graph containing (exactly) one $4$-cycle.

\medskip \centerline{*}

Let us finish this section with a comparison of parameter choices.
In classic graph coloring theory, allowing an additional color results in a weaker coloring, or, with a more positive attitude, allows colorings of a bigger class of graphs.
In the threshold coloring setting additional colors might not allow coloring of additional graphs unless we also allow threshold to grow. 

This motivates us to define a relation on threshold parameter pairs by writing
$(r_1,t_1) \leq (r_2,t_2)$ if and only if there exists an increasing injective mapping 
$\varphi: \{0,\ldots,r_1-1\} \rightarrow \{0,\ldots,r_2-1\}$
so that for every par of integers $a,b \in \{0,\ldots,r_1-1\}$ we have 
$|a - b| \le t_1$ if and only if $|\varphi(a) - \varphi(b)| \le t_2$.

Clearly, $\leq$ is a partial order. On the other hand, it is not a total order, as there are incomparable parameter pairs; take,
for example, $(11,1)$ and $(18,4)$. Then $(18,4) \not\leq (11,1)$ as there is no injective mapping from the first $18$ to the first $11$ nonnegative integers.
Assuming, to the contrary, that there was a mapping $\varphi$ as above certifying $(11,1) \leq (18,4)$, we would know that
$\varphi(0)$, $\varphi(2)$, $\varphi(4)$, $\varphi(6)$, $\varphi(8)$, $\varphi(10)$ would be an increasing sequence where consecutive members differ by at least $5$;
this implies $\varphi(10) \geq 25$, contradiction.

It is easy to see that if $(r_1,t_1) \leq (r_2,t_2)$ and $G$ admits an $(r_1,t_1)$-threshold-coloring $c$ with respect to some near-far-labeling $(N,F)$ then
$\varphi \circ c$ is an $(r_2,t_2)$-threshold-coloring with respect to $(N,F)$. 
By taking $\varphi$ to be the identity we see immediately that $(r_1,t) \leq (r_2,t)$ if $r_1 \leq r_2$,
and by defining $\varphi(a):=\lambda a$ we see that $(r,t) \leq (\lambda r, \lambda t)$ for each positive integer $\lambda$.

The important fact is that $(r_1,t_1)$ and $(r_2,t_2)$ always have a common upper bound with respect to $\leq$:
Take $t$ to be the least common multiple of $t_1$ and $t_2$
then $\lambda_1:=t/t_1$ and $\lambda_2:=t/t_2$ are integers and setting $r:=\max\{\lambda_1 r_1,\lambda_2 r_2\}$ we see that
$(r_1,t_1) \leq (\lambda_1 r_1,\lambda_1 t_1)=(\lambda_1 r_1,t) \leq (r,t)$ and, analogously, $(r_2,t_2) \leq (r,t)$.
The next lemma states that we can always find a common upper bound $(r,t)$ with respect to $\leq$ such that $t$ is the maximum of $t_1,t_2$.
Results of this type allow to split an argument showing that there exist $r,t$ such that there is an $(r,t)$-threshold-coloring with respect to every near-far-labeling $(N,F)$ into
cases according to structural properties of $(N,F)$: If for each labeling $\pi:=(N,F)$ we find an $(r_\pi,t_\pi)$-coloring of $G$
then just take $(r,t)$ as a common upper bound of all the $(r_\pi,t_\pi)$ --- and $G$ will be 
$(r,t)$-total-threshold-colorable. 

\begin{lem}
  Let $(r_1,t_1)$ and $(r_2,t_2)$ be parameter pairs, where $t_1 \leq t_2$.
  Then $(r_1,t_1) \le (r,t_2)$ and $(r_2,t_2) \le (r,t_2)$ for
  \[ r = \max \{r_2, (r_1 \ddiv (t_1 +1)) \cdot 
   (t_2+1) + (r_1 \bmod (t_1+1)), \]
  where $\ddiv$ and $\bmod$ denote the standard integer quotient and remainder, respectively.
\end{lem}

{\bf Proof.}
Clearly $(r_2,t_2) \le (r,t_2)$, as $r_2 \le r$. 

The mapping 
$$\varphi: x \mapsto (x \ddiv (t_1 +1)) \cdot (t_2 + 1) + (x \bmod (t_1+1))$$
is an increasing injective mapping from $\{0,\ldots,r_1-1\}$ to
\begin{eqnarray*}
& & \{\varphi(0),\ldots, \varphi(r_1-1)\} \\
& \subseteq & \{0,\ldots, ((r_1-1) \ddiv (t_1 +1)) \cdot (t_2+1) + ((r_1-1) \bmod (t_1+1)) \} \\
& \subseteq & \{0,\ldots, (r_1 \ddiv (t_1 +1)) \cdot (t_2+1) + (r_1 \bmod (t_1+1)) - 1 \} \\
& \subseteq & \{0,\ldots,r-1\}.
\end{eqnarray*}

For $a,b \in \{0,\ldots,r_1-1\}$ which satisfy $|a-b|=t_1+1$ we also have 
$|\varphi(a)-\varphi(b)| = t_2+1$. As $\varphi$ is strictly increasing, we also have
$(r_1,t_1) \le (r,t_2)$. 
\hspace*{\fill}$\Box$

To continue the above example we see that, at the expense of increasing the range, we can ``interpret'' both an $(11,1)$- and $(18,4)$-threshold-coloring as a $(26,4)$-threshold-coloring
since $(11 \ddiv (1+1))\cdot (4+1) + (11 \bmod (1+1))=26$, whereas the more elementary common upper bound would be $(44,4)$.
Observe that $(5,1) \leq (11,4) \leq (13,4)$, as (for the first inequality) $(5 \ddiv (1+1) \cdot (4+1)+(5 \bmod (1+1))=11$, so that 
the two parameter pairs in Lemma \ref{zigzag} are comparable, too.

\section{Ladders and M\"obius ladders}
\label{SLadders}

\begin{figure}
\begin{center}
\begin{tikzpicture}[scale=1.45]
\labeledladder{3} 
\end{tikzpicture}
\hspace{5mm}
\begin{tikzpicture}[scale=1.45]
\labeledprism{4}
\end{tikzpicture}
\end{center}
\caption{The Ladder $L_4$ and the prism $C_5 \cpg K_2$, shown with notation of edges, vertices, and squares.}
\label{ladderprismlabels}
\end{figure}
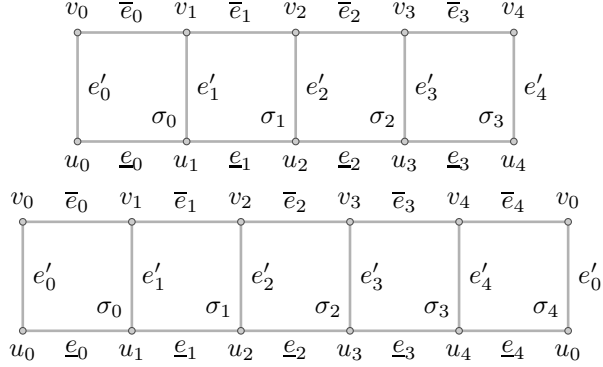

Recall that the ladder $L_n$, $n \geq 0$, is defined to be the cartesian product $P_n \cpg K_2$ of the path $P_n$ of length $n$ and
the complete graph $K_2$ on two vertices. The cartesian product $C_n \cpg K_2$ of a cycle $C_n$ of length $n \geq 3$ and $K_2$
is called a {\em prism over $C_n$}. An edge in $L_n$ or $C_n \cpg K_2$ is called a {\em spoke} if it is an edge of a $K_2$-fibre,
the remaining edges are called {\em peripheral edges}.
A {\em square} $\sigma$ in $L_n$ or $C_n \cpg K_2$ is a $4$-cycle containing two spokes and two peripheral edges.
The notational conventions we follow when studying edges, vertices, and squares in ladders or prisms are depicted in Figure \ref{ladderprismlabels};
accordingly, for a square $\sigma=\sigma_i$ we get $V(\sigma)=\{v_i,v_{i+1},u_i,u_{i+1}\}$ and $E(\sigma)=\{\overline{e}_i,\underline{e}_i,e'_i,e'_{i+1}\}$.

Note that the prism over $C_n$ can be obtained from the ladder $L_n$ by identifying a pair of extremal spokes (the ones contained in only one square),
so that the peripheral edges span a pair of {\em peripheral $n$-cycles}. 
The other possible identification of boundary spokes results in a {\sc M\"obius} ladder $M_n$.
We define spokes, peripheral edges, and squares in $M_n$ in an analogous way. We shall however denote the vertices of $M_n$ by 
$v_0, \ldots, v_{2n-1}$ along the unique peripheral cycle.

\begin{lem}
\label{ladder}
$L_n$ is $(5,1)$-total-threshold-colorable.
\end{lem}

{\bf Proof.}
Let us set $A:=\{v_0,v_4,v_8,\ldots,u_2,u_6,u_{10},\ldots\}$, 
and, naturally, $B:= V(L_n) \setminus A$. Let $(N,F)$ be an arbitrary edge labeling.
We shall show that the conditions of Lemma~\ref{zigzag} are satisfied:
As $A$ is $2$-independent, (i) holds irrespective of the labeling. Similarly, $B$ induces a path, hence (ii) is also satisfied.
Finally, every vertex $v \in B$ has a unique neighbor in $A$, hence (iv) (and, trivially, (iii)) hold with $M=\emptyset$.
Lemma~\ref{zigzag} implies that $L_n$ is $(5,1)$-total-threshold-colorable with respect to the chosen labeling. As the labeling was arbitrary the proof is finished.
\hspace*{\fill}$\Box$

Next we shall prove that the {\sc M\"obius} ladders are not total threshold colorable.

{\bf Proof of Theorem \ref{moebiusmain}.}


Although the case $K_{3,3}=M_3$ follows from the general result, let us give a separate argument. 
Let us denote the vertices $K_{3,3}$ by $v_0,v_1,v_2,v_3,v_4,v_5$, where two vertices are adjacent if and only if the parity of their indices differ.
Set $F=\{v_0v_3, v_1v_4, v_2v_5\}$, the remaining edges being labeled near (so $N:=E(K_{3,3}) \setminus F$), and assume that $c$ is a threshold coloring with respect to $(N,F)$.
As $v_1 v_0 \in N$ and $v_1 v_4 \in F$ 
we have $c(v_0) \ne c(v_4)$. By a similar argument we infer that vertices $v_0, v_2, v_4$ receive three different colors, and by symmetry we may assume $c(v_0) < c(v_2) < c(v_4)$. Now the color $c(v_5)$ cannot be both far from the middle color $c(v_2)$ and close to both extremal colors $c(v_0)$ and $c(v_4)$. 

Let us finish with a general argument.
Let $(N,F)$ be the near-far-labeling of $M_n$ in which exactly the spokes are far edges. 
Assume that $c$ is a threshold coloring with respect to $(N,F)$, for an appropriate choice of parameters $r,t$.
Without loss of generality we may assume that $c(v_0) > c(v_{n})$.
By Lemma~\ref{init} we infer that $c(v_1) > c(v_{n+1})$,
and inductively also $c(v_{k}) > c(v_{n+k})$ for every integer $k$,
where the addition in indices is taken modulo $2n$. By setting $k=n$ we have $c(v_n) > c(v_0)$ which is absurd. 
\hspace*{\fill}$\Box$

Let us use Theorem \ref{moebiusmain} and Lemma \ref{ladder} to proof that being total threshold colorable cannot be characterize in terms of a finite set forbidden subgraphs.

{\bf Proof of Corollary \ref{cormoebius}.}

Assume, to the contrary, that there does exist a \emph{finite} set $\mathfrak F$ such that
a graph is total threshold colorable if and only if it contains no [induced] subgraph isomorphic to one of those in $\mathfrak F$.
Then every single graph from $\mathfrak F$ is not total threshold colorable. Let $n:=\max\{|V(G)|:\, G \in {\mathfrak F}\}$,
and consider the M\"obius lader $M_{n+1}$, which has $2n+2$ vertices.
By assumption, $M_{n+1}$ contains a subgraph $H$ isomorphic to a graph from $\mathfrak F$.
Then $H$ is not total threshold colorable and has at most $n$ vertices.
By the pidgeon hole principle, there exists a spoke $v_k v_{k+(n+1)}$ of $M_{n+1}$ such that neither $v_k$ nor $v_{k+(n+1)}$ are in $V(H)$,
that is, $H$ is a subgraph of the graph $M_{n+1}-\{v_k,v_{k+(n+1)}\}$, which is isomorphic to $L_n$.
By Lemma \ref{ladder}, $L_n$, and, hence, $H$, is total threshold colorable, contradiction.
\hspace*{\fill}$\Box$

\section{Prisms}
\label{SPrisms}

Before turning attention to the proof of Theorem \ref{prismain}, let us give a handful of definitions.
Consider a square $\sigma$ in either $L_n$ or $C_n \cpg K_2$, that is, a $4$-cycle containing exactly two consecutive spokes.
A square is {\em balanced} if it contains exactly two far and two near edges, and is {\em unbalanced} otherwise.
An edge $e$ of an unbalanced square is a {\em deviator} if the remaining three edges are labeled different from $e$.
Note that a deviator exists only if and only if the edges in a square split $3:1$ according to their labels. 

\begin{figure}
\begin{center}
\begin{tikzpicture}[scale=0.70]
\ladder{1};
\draw[near] (0,0) -- +(1,0);
\draw[near] (0,0) -- +(0,1);
\draw[near] (1,1) -- +(-1,0);
\draw[near] (1,1) -- +(0,-1);
\end{tikzpicture}
\hspace{15mm}
\begin{tikzpicture}[scale=0.70]
\ladder{1};
\draw[near] (0,0) -- +(1,0);
\draw[near] (0,0) -- +(0,1);
\draw[near] (1,1) -- +(-1,0);
\draw[far] (1,1) -- +(0,-1);
\end{tikzpicture}
\hspace{15mm}
\begin{tikzpicture}[scale=0.70]
\ladder{1};
\draw[near] (0,0) -- +(1,0);
\draw[far] (0,0) -- +(0,1);
\draw[near] (1,1) -- +(-1,0);
\draw[far] (1,1) -- +(0,-1);
\end{tikzpicture}
\hspace{15mm}
\begin{tikzpicture}[scale=0.70]
\ladder{1};
\draw[near] (0,0) -- +(1,0);
\draw[near] (0,0) -- +(0,1);
\draw[far] (1,1) -- +(-1,0);
\draw[far] (1,1) -- +(0,-1);
\end{tikzpicture}
\end{center}
\caption{Examples of unbalanced even, unbalanced odd, (balanced) parallel, and (balanced) nonparallel squares. Near edges are thick, far edges are dashed.}
\label{fig:square:types}   
\end{figure}
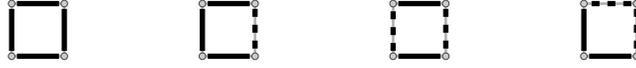

A square is {\em even} if it contains an even number of both near and far edges, and {\em odd} if its edges split $3:1$ with respect to their labels. 
A balanced square is {\em parallel} if the two near edges form a matching (also the far edges form a matching) and {\em nonparallel} if its near edges are consecutive.
We shall adopt the convention to draw near edges in thick black and far edges dashed, examples are shown in Figure~\ref{fig:square:types}.

{\bf Proof of Theorem~\ref{prismain}}

Let $G:=C_n \cpg K_2$, $n \ge 3$, be a prism over $C_n$, equipped with an arbitrary edge labeling $\pi=(N,F)$.
We have to prove that $G$ admits $(26,4)$-threshold-coloring with respect to $\pi$.

 

To this end, we define some special sets of far edges.
For $i,j \in \{0,\dots,n-1\}$, a set $S \subseteq F$ is called a {\em half-cut 
across the squares $\sigma_i,\sigma_{i+1},\dots,\sigma_k$} or simply a {\em half-cut},
if it is equal to either $\{\underline{e}_i,\overline{e_k},e_{i+1}',e_{i+2}',\dots,e_k'\}$ or $\{\overline{e}_i,\underline{e_k},e_{i+1}',e_{i+2}',\dots,e_k'\}$,
where the indices are taken modulo $n$; if $i=k$ then the half-cut across $\sigma_i$ consists of $\{\underline{e}_i,\overline{e}_i\}$.
Observe that $G-S$ is always isomorphic to a subgraph of some ladder $L_{n'}$ with $n' \geq n$.
A set $S \subseteq F$ is {\em useful} if there exists a partition $\{A,B\}$ of $V(G)$ such that each of $G[A],G[B]$ is isomorphic to a (not necessarily connected) subgraph of a ladder
and an edge of $G$ is in $S$ if and only if it has endvertices in both $A$ and $B$. In particular, $S$ is a (not necessarily minimal) cut, and we sometimes call it a {\em useful cut}.

{\bf Claim 1.} The union of two disjoint half-cuts is a useful cut.

Let us first look at an illustrative example, where the half-cuts 
$S$ and $T$ each contain two edges. Without loss of generality we may assume $S=\{\underline{e}_i, \overline{e}_i\}$ and $T=\{\underline{e}_j, \overline{e}_j\}$, for some choice of $i < j$. In this case $G-S-T$ is a union of two disjoint ladders of lengths $j-i-1$ and $n-j+i-1$.

In what follows we may without loss of generality assume that
$|S| \ge |T|$, $|S| \ge 3$, $S$ is a half-cut across $\sigma_{n-1}, \sigma_0,\ldots, \sigma_i$, and also $\{\overline{e}_{n-1}, \underline{e}_i\} \subseteq S$. 

If $T$ contains two peripheral edges from a square $\sigma_j$, where $j < i$, then $T$ contains no spokes, as $S$ contains both $e'_j$ and $e'_{j+1}$. In this case we can replace half-cuts $S$ and $T$ with a pair of alternative half-cuts across $\sigma_{n-1},\ldots, \sigma_j$ and $\sigma_j,\ldots,\sigma_i$, respectively, whose union is equal to $S \cup T$.

Therefore we can also assume that $T$ is a half-cut across squares $\sigma_j,\ldots, \sigma_k$, where $j \ge i$ and $k \le n-1$.

Now if $\overline{e}_j \in T$ (and consequently $\underline{e}_k \in T$) we set $A=\{v_0,\ldots,v_j, u_{i+1},\ldots,u_{k+1}\}$, and in case
$\underline{e}_j \in T$ (and  $\overline{e}_k \in T$)
we set $A=\{v_0,\ldots,v_{k+1}, u_{i+1},\ldots,u_{j}\}$. Note also that in the latter case $j>i$, as $\underline{e}_i \not\in T$. 
We also let $B = V(G) \setminus A$.

Now $A$ and $B$ each induce a pair of paths in peripheral cycles of $G$
and every spoke $e'_\ell \not \in S \cup T$ has both endvertices either in $A$ (if $i< \ell \le j$) or $B$ (otherwise, $k < \ell \le n-1$). This implies that both $G[A]$ and $G[B]$ are isomorphic to subgraphs of ladders, which in turn proves Claim 1.

{\bf Claim 2.} If there exists a useful cut in $G$ then $G$ is $(11,1)$-threshold-colorable with respect to $(N,F)$.

Let $S \subseteq F$ be useful and $\{A,B\}$ a partition of $V(G)$ such that $G[A],G[B]$ are subgraphs of a ladder and
an edge of $G$ is in $S$ if and only if it has endvertices in both $A$ and $B$.
By Lemma \ref{ladder} there exists a $(5,1)$-threshold-coloring $a$ of $G[A]$ with respect to $(N \cap E(G[A]),F \cap E(G[A]))$,
and a $(5,1)$-threshold-coloring $b$ of $G[B]$ with respect to $(N \cap E(G[B]),F \cap E(G[B]))$.
We may assume that $a$ uses colors $\{0,\dots,4\}$ and $b$ uses colors from $\{6,\dots,10\}$.
Hence $a \cup b$ is an $(11,1)$-threshold-coloring of $G$ with respect to $(N,F)$.
This proves Claim 2.

{\bf Claim 3.} Suppose that $G$ and $(N,F)$ have the following properties.
\begin{itemize}
  \item[(i)] There are no two disjoint half-cuts.
  \item[(ii)] Every vertex is incident with a near edge.
  \item[(iii)] $G$ admits no parallel square with far spokes unless $n=3$.
\end{itemize}
Then $G$ is $(13,4)$-threshold-colorable with respect to $(N,F)$.

We start the proof of Claim 3 by showing the following.

{\bf Subclaim.} $G$ admits an unbalanced square unless $n=3$.



Let $n \geq 4$ and assume, to the contrary, that every square of $G$ is balanced.
By (iii) we may assume that every square is either parallel with near spokes or nonparallel. In particular, no two consecutive spokes are far.

Let $e_0$ be a far spoke contained in consecutive squares $\sigma$ and $\sigma'$.
As $\sigma$ and $\sigma'$ are nonparallel, $e_0$ is adjacent to a pair of far edges $e$, $e'$ contained in $\sigma$ and $\sigma'$, respectively.
By (ii) the edges $e_0, e,e'$ form a half cut, as they are not incident with a common vertex.
A pair of nonconsecutive far spokes thus gives rise to a pair of edge-disjoint half-cuts, which is by (i) not possible. 

Hence there exists a parallel square $\sigma_0$ with near spokes, which by itself contains a half-cut. 
Let $\sigma_0, \sigma_1, \sigma_2$ be a sequence of three consecutive squares.
As $\sigma_1$ is nonparallel 
(otherwise we have a pair of disjoint half cuts), 
the spoke $e_0 \in E(\sigma_1) \cap E(\sigma_2)$ is a far spoke.
By above argument it is contained in 
another 
half-cut, which is again absurd. 
This proves our subclaim.


In what remains we shall split our analysis according to $n \bmod 4$.

{\bf Case 1.} $n=4k+3$ for some integer $k \geq 0$. 
 
Assume first that $G$ contains an unbalanced square $\sigma=v_0 u_0 u_{4k+2} v_{4k+2}$ (as it does in the case $n \geq 7$ by our subclaim).
Let us define
$A:=\{v_0,v_4$, $\ldots$, $v_{4k}$, $u_2$, $u_6$, $\ldots$, $u_{4k+2}\}$, and let $B:= V(G) \setminus A$. 
We shall see that we can apply Lemma~\ref{zigzag} in this case.

Let us first determine the set of edges of $M$ satisfying (iii) of Lemma \ref{zigzag}.
As $\sigma$ is unbalanced there may exist a deviator $e$ in $\sigma$.
In this case we set $M:=\{e\}$, and in the remaining case let $M:=\emptyset$.
Clearly, $A$ is an independent set, hence (i) of Lemma \ref{zigzag} holds.

\begin{figure}
\begin{center}
\begin{tikzpicture}[scale=0.85]
\unbalancedsquare{4};
\ladder{7}
\draw[zigzag] (0,0) -- +(1,0);
\draw[zigzag] (1,0) -- +(1,0);
\draw[zigzag] (2,1) -- +(1,0);
\draw[zigzag] (4,0) -- +(1,0);
\draw[zigzag] (5,1) -- +(1,0);
\draw[zigzag] (6,1) -- +(1,0);

\draw[zigzag] (0,0) -- +(0,1);
\draw[zigzag] (2,0) -- +(0,1);
\draw[zigzag] (5,0) -- +(0,1);
\draw[zigzag] (7,0) -- +(0,1);


\node at (3.65, 0.35) {$\sigma_0$};
\node[blackvertex] at (1,1) {};
\node[blackvertex] at (3,0) {};
\node[blackvertex] at (4,1) {};
\node[blackvertex] at (6,0) {};

\end{tikzpicture}
\hspace{10mm}
\begin{tikzpicture}[scale=0.85]
\ladder{7}
\draw[zigzag] (0,0) -- +(1,0);
\draw[zigzag] (1,0) -- +(1,0);
\draw[zigzag] (2,1) -- +(1,0);
\draw[zigzag] (4,0) -- +(1,0);
\draw[zigzag] (5,1) -- +(1,0);
\draw[zigzag] (6,1) -- +(1,0);
\draw[zigzag] (0,0) -- +(0,1);
\draw[zigzag] (2,0) -- +(0,1);
\draw[zigzag] (5,0) -- +(0,1);
\draw[zigzag] (7,0) -- +(0,1);
\draw[far] (3,0) -- +(0,1);
\draw[far] (3,1) -- +(1,0);
\draw[near] (4,0) -- +(0,1);
\draw[near] (3,0) -- +(1,0);
\node[blackvertex] at (1,1) {};
\node[blackvertex] at (3,0) {};
\node[blackvertex] at (4,1) {};
\node[blackvertex] at (6,0) {};

\node at (3.65, 0.35) {$\sigma_0$};
\end{tikzpicture}
\end{center}
\caption{An unbalanced square and a balanced nonparallel square $\sigma$ in $C_7 \cpg K_2$,
and square-wave paths. Black vertices will receive color $0$, in an unbalanced square before possible perturbation.}
\label{fig:mod3}
\end{figure}
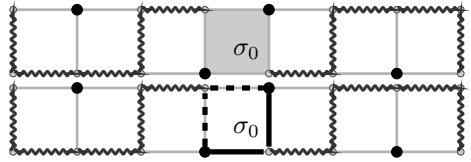

The set $B$ induces a path $P$, containing every second spoke, see Figure~\ref{fig:mod3}.
We shall call such a path a {\em zigzag} path, and it trivially satisfies (ii) of Lemma \ref{zigzag}.
As $G$ is cubic, every interior vertex of $P$ is adjacent to exactly one vertex of $A$.
It is the endvertices of $P$ that might interfere with (iv). However this is not the case if $M= \emptyset$ as then every edge of $\sigma$ has the same label.
But if $M=\{e\}$ then (iv) is also satisfied, as $e$ does not enter (iv) of Lemma \ref{zigzag}.
By Lemma~\ref{zigzag} there exists a $(13,4)$-threshold-coloring.

We are now in the case that all squares are balanced.
Although it would suffice to restrict to the case $n=3$ from here on, we give a more general argument.
Note that we can apply Lemma~\ref{zigzag}, with $M=\emptyset$, also in the case where
$G$ contains a balanced nonparallel square $\sigma$, see Figure~\ref{fig:mod3}.
In this case we can achieve that $V(\sigma) \cap B$ consists of vertices incident with edges of the same label along $\sigma$
(otherwise we swap the roles of $u$ and $v$ in the definition of $A$).

Hence we may assume that every square of $G$ is balanced and parallel.
Consequently, either every spoke is a near edge or every spoke is a far edge.
The former option is by (i) not possible, as every square balanced parallel square with near spokes admits a half-cut across it.
The latter option allows a $(2,0)$-threshold-coloring of $G$ by colors $\{0,1\}$ in such a way
that every near edge is incident with vertices of the same color. 
As $(2,0) \leq (13,4)$, this proves Claim 3 in Case 1. 

Before drilling into the remaining cases let us focus on a common feature. 
Given an unbalanced or balanced nonparallel square $\sigma$ in $G$ with $n=4k+3$, $k \geq 0$, 
we were able to connect a square-wave path $P$ to two antipodal vertices of $\sigma$, so that the vertices of $V(G) \setminus V(P)$ formed an independent set.
In case $n \bmod 4 \not= 3$ we need an additional tool to succeed, as the square-wave path $P$ cannot be properly joined to antipodal vertices of $\sigma$.
In case $n \bmod 4 = 0$ or $1$ we will have to change the phase of $P$ so that the proper connection to $\sigma$ can be accomplished.
In case $n \bmod 4 =2$ we shall try to use two disjoint square-wave paths, connected to a pair of squares.  

\begin{figure}
\begin{center}
\begin{tikzpicture}[scale=0.85]
\ladder{12}
\draw[zigzag] (0,0) -- +(1,0);
\draw[zigzag] (1,1) -- +(1,0);
\draw[zigzag] (2,1) -- +(1,0);
\draw[zigzag] (3,0) -- +(1,0);
\draw[zigzag] (4,0) -- +(1,0);
\draw[zigzag] (5,1) -- +(1,0);
\draw[zigzag] (6,1) -- +(1,0);
\draw[zigzag] (7,0) -- +(1,0);
\draw[zigzag] (8,0) -- +(1,0);
\draw[zigzag] (9,1) -- +(1,0);
\draw[zigzag] (10,1) -- +(1,0);
\draw[zigzag] (11,0) -- +(1,0);
\draw[zigzag] (1,0) -- +(0,1);
\draw[zigzag] (3,0) -- +(0,1);
\draw[zigzag] (5,0) -- +(0,1);
\draw[zigzag] (7,0) -- +(0,1);
\draw[zigzag] (9,0) -- +(0,1);
\draw[zigzag] (11,0) -- +(0,1);
\node[blackvertex] at (0,1) {};
\node[blackvertex] at (4,1) {};
\node[blackvertex] at (8,1) {};
\node[blackvertex] at (12,1) {};
\node[blackvertex] at (2,0) {};
\node[blackvertex] at (6,0) {};
\node[blackvertex] at (10,0) {};
\end{tikzpicture}
\vspace{5mm}

\begin{tikzpicture}[scale=0.85]
\ladder{12}
\draw[near] (0,1) -- +(1,0);
\draw[near] (3,0) -- +(1,0);
\draw[near] (4,0) -- +(1,0);
\draw[zigzag] (0,0) -- +(1,0);
\draw[zigzag] (1,0) -- +(1,0);
\draw[zigzag] (2,1) -- +(1,0);
\draw[zigzag] (3,1) -- +(1,0);
\draw[zigzag] (4,1) -- +(1,0);
\draw[zigzag] (5,1) -- +(1,0);
\draw[zigzag] (6,0) -- +(1,0);
\draw[zigzag] (7,0) -- +(1,0);
\draw[zigzag] (8,0) -- +(1,0);
\draw[zigzag] (8,1) -- +(1,0);
\draw[zigzag] (9,1) -- +(1,0);
\draw[zigzag] (10,1) -- +(1,0);
\draw[zigzag] (2,0) -- +(0,1);
\draw[zigzag] (6,0) -- +(0,1);
\draw[zigzag] (8,0) -- +(0,1);
\draw[zigzag] (9,0) -- +(0,1);
\draw[zigzag] (11,0) -- +(0,1);
\draw[zigzag] (11,0) -- +(1,0);
\node at (8.65, 0.35) {$\sigma$};
\node[blackvertex] at (0,1) {};
\node[blackvertex] at (1,1) {};
\node[blackvertex] at (7,1) {};
\node[blackvertex] at (12,1) {};
\node[blackvertex] at (3,0) {};
\node[blackvertex] at (4,0) {};
\node[blackvertex] at (5,0) {};
\node[blackvertex] at (10,0) {};
\end{tikzpicture}
\end{center}
\caption{Square-wave path (top) and phase delays using near edges and even squares (bottom). Every near edge (shown thick) delays the wave by 1, so does an even square $\sigma$. Again black vertices indicate color $0$.}
\label{fig:delay}
\end{figure}
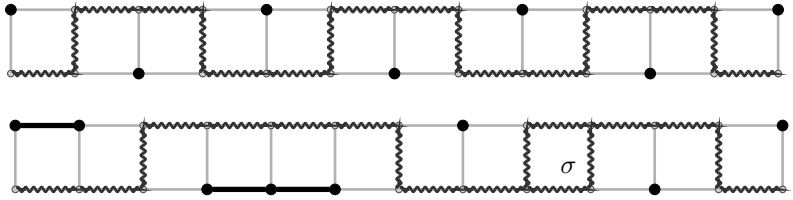

There are two different tools to delay the square-wave path, see Figure~\ref{fig:delay}. 
A correctly placed near edge $e$ away from $P$ may be used to allow a pair of adjacent vertices in $A$.
Similarly, an even square may be used to beef up the run of $P$, as we allow vertices of $B$ to induce even cycles.

Note also that, in the remaining cases, $n \geq 4$, which implies that no balanced square has both spokes labeled far by (iii),
and that there always exists an unbalanced square by our subclaim.

{\bf Case 2.} $n=4k+4$ for some integer $k \geq 0$.

By our subclaim, $G$ admits an unbalanced square. 
Without loss of generality, we may assume $\sigma_0=v_0 v_1 u_1 v_0$ to be unbalanced.
By (i), not every peripheral edge incident with $\sigma_0$ is a far edge,
hence, by symmetry, we can assume that $v_1v_2 \in N$.

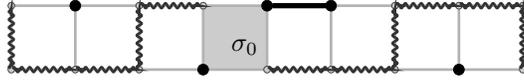
\begin{figure}
\begin{center}
\begin{tikzpicture}[scale=0.85]
\unbalancedsquare{4};
\ladder{8}
\draw[zigzag] (0,0) -- +(1,0);
\draw[zigzag] (1,0) -- +(1,0);
\draw[zigzag] (2,1) -- +(1,0);
\draw[zigzag] (4,0) -- +(1,0);
\draw[near] (4,1) -- +(1,0);
\draw[zigzag] (5,0) -- +(1,0);
\draw[zigzag] (6,1) -- +(1,0);
\draw[zigzag] (7,1) -- +(1,0);
 
\draw[zigzag] (0,0) -- +(0,1);
\draw[zigzag] (2,0) -- +(0,1);
\draw[zigzag] (6,0) -- +(0,1);
\draw[zigzag] (8,0) -- +(0,1);

\node[blackvertex] at (1,1) {};
\node[blackvertex] at (3,0) {};
\node[blackvertex] at (4,1) {};
\node[blackvertex] at (5,1) {};
\node[blackvertex] at (7,0) {};

\node at (3.65, 0.35) {$\sigma_0$};
\end{tikzpicture}
\end{center}
\caption{Unbalanced square $\sigma_0$ in $C_8 \cpg K_2$.}
\label{fig:mod4}
\end{figure}

Let us define $A:=\{v_1, v_2, v_6,\ldots,v_{4k+2},u_4,u_8,\ldots,u_{4k+4}=u_0\}$, and let $B:= V(G) \setminus A$, see the illustration on Figure~\ref{fig:mod4}.

It is possible that there exists a deviator $e$ of $\sigma_0$;
in that case, we set $M:=\{e\}$, and otherwise $M:=\emptyset$. 
We claim that $G,(N,F),A,B,M$ satisfy the conditions of Lemma~\ref{zigzag}.
Clearly $G[A]$ contains a single edge $v_1v_2 \in N$, hence (i) of Lemma \ref{zigzag} holds.
Next $G[B]$ is a path, so (ii) of Lemma \ref{zigzag} is trivially satisfied, and so is (iii) as $M$ contains at most one edge.
Finally, (iv) of Lemma \ref{zigzag} follows as every internal vertex of $G[B]$ has a single neighbor in $A$,
and endvertices lie on an unbalanced cycle and a possibly problematic edge $e$ is in $M$.
Hence $G$ admits a $(13,4)$-threshold-coloring with respect to $(N,F)$.
This settles Claim 3 in Case 2.

{\bf Case 3.} $n=4k+6$ for some integer $k \geq 0$.

Again, we may without loss of generality assume that $\sigma_0$ is an unbalanced square.

Assume first that $\sigma_3$ is either an unbalanced square or a nonparallel balanced square
(in the latter case we may by symmetry assume that edges $e'_3$ and $\underline{e_3}$ have the same label). 
Let $f_0$ be the deviator of $\sigma_0$ (if it exists), and $f_3$ be the deviator of $\sigma_3$ (if it exists).
Let $M$ be the set of deviators of squares $\sigma_0$ and $\sigma_3$ (so $|M| \in \{0,1,2\}$).

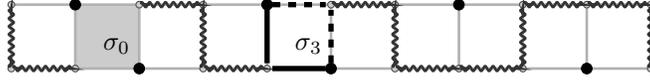
\begin{figure}
\begin{center}
\begin{tikzpicture}[scale=0.85]
\unbalancedsquare{2};
\ladder{10}
\draw[zigzag] (0,0) -- +(1,0);
\draw[zigzag] (2,1) -- +(1,0);
\draw[zigzag] (3,0) -- +(1,0);
\draw[zigzag] (5,1) -- +(1,0);
\draw[zigzag] (6,0) -- +(1,0);
\draw[zigzag] (7,0) -- +(1,0);
\draw[zigzag] (8,1) -- +(1,0);
\draw[zigzag] (9,1) -- +(1,0);
\draw[zigzag] (0,0) -- +(0,1);
\draw[zigzag] (3,0) -- +(0,1);
\draw[zigzag] (6,0) -- +(0,1);
\draw[zigzag] (8,0) -- +(0,1);
\draw[zigzag] (10,0) -- +(0,1);
\draw[near] (4,0) -- +(1,0);
\draw[near] (4,0) -- +(0,1);
\draw[far] (4,1) -- +(1,0);
\draw[far] (5,0) -- +(0,1);
\node at (1.65, 0.35) {$\sigma_0$};
\node at (4.65, 0.35) {$\sigma_3$};
\node[blackvertex] at (1,1) {};
\node[blackvertex] at (2,0) {};
\node[blackvertex] at (4,1) {};
\node[blackvertex] at (5,0) {};
\node[blackvertex] at (7,1) {};
\node[blackvertex] at (9,0) {};

\end{tikzpicture}
\end{center}
\caption{Unbalanced square $\sigma_0$ and a balanced nonparallel square $\sigma_3$ in $C_{10} \cpg K_2$.}
\label{fig:mod4:2}
\end{figure}

Let $A= \{u_1, v_3, u_4 , u_8, \ldots, u_{4k+4}, v_6, v_{10}, \ldots, v_{4k+6} = v_0 \}$ and  $B = V(G) \setminus A$.
Now $G$, $(N,F)$, $A$, $B$, $M$ satisfy the conditions of Lemma~\ref{zigzag},
as the edges of $M$ form an induced matching.
Hence there exists a $(13,4)$-threshold-coloring of $G$ with respect to $(N,F)$, see Figure~\ref{fig:mod4:2}.

Hence we may assume that  $\sigma_3$ is a balanced parallel square.
By (iii) its peripheral edges are labeled far, and they thus form a half-cut across $\sigma_3$.
Another unbalanced square $\sigma_i$, $i \not= 0$, implies 
that the peripheral edges of $\sigma_{i+3}$ would form a half-cut disjoint from the first; this is impossible by (i).
Therefore $\sigma_0$ is the {\em only} unbalanced square, $\sigma_3$ is parallel and admits a half-cut across it, and all the remaining squares are balanced nonparallel. 

As $\sigma_2$ is balanced nonparallel and 
$e_3'$ 
is near, the spoke $e'_2$ is a far edge.
Together with the remaining two far edges from $E(\sigma_1) \cup E(\sigma_2)$ we either have a vertex incident with three far edges or another (disjoint from the one in $\sigma_3$) half-cut.
By (i) and (ii) neither is possible and Claim 3 is settled in Case 3.

{\bf Case 4.} $n=4k+5$ for some integer $k \geq 0$.

This is our final case, and also the most difficult one.
While it will be rather easy to take care of the longer prisms, the shortest prism in this case, that is, $G=C_5 \cpg K_2$, will prove to be a difficult beast. 

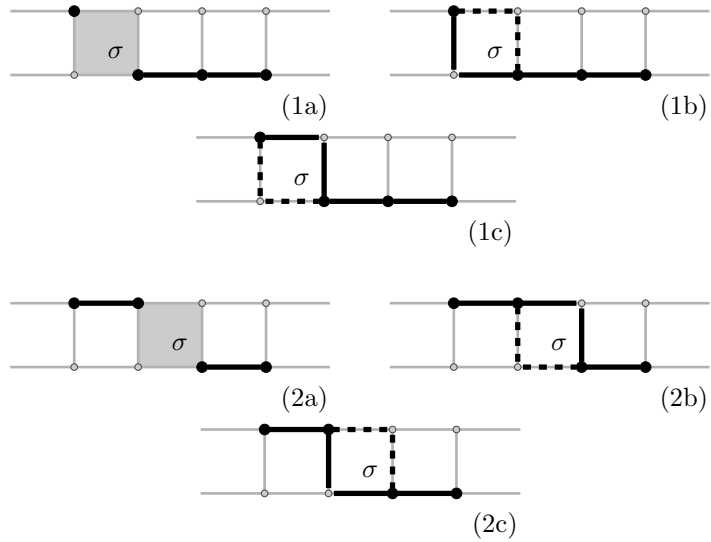
\begin{figure}
\begin{center}
\begin{tikzpicture}[scale=0.85]
\unbalancedsquare{1};
\ladderplus{3}
\draw[near] (1,0) -- +(1,0);
\draw[near] (2,0) -- +(1,0);
\node[blackvertex] at (0,1) {};
\node[blackvertex] at (1,0) {};
\node[blackvertex] at (2,0) {};
\node[blackvertex] at (3,0) {};
\node at (0.65, 0.35) {$\sigma$};
\node at (3.6,-0.5) {(1a)};
\end{tikzpicture}
\hspace{13pt}
%
\begin{tikzpicture}[scale=0.85]
\ladderplus{3}
\draw[near] (0,0) -- +(1,0);
\draw[near] (0,0) -- +(0,1);
\draw[far] (1,1) -- +(-1,0);
\draw[far] (1,1) -- +(0,-1);
\draw[near] (1,0) -- +(1,0);
\draw[near] (2,0) -- +(1,0);
\node[blackvertex] at (0,1) {};
\node[blackvertex] at (1,0) {};
\node[blackvertex] at (2,0) {};
\node[blackvertex] at (3,0) {};
\node at (0.65, 0.35) {$\sigma$};
\node at (3.6,-0.5) {(1b)};
\end{tikzpicture}
\hspace{13pt}
%
\begin{tikzpicture}[scale=0.85]
\ladderplus{3}
\draw[far] (0,0) -- +(1,0);
\draw[far] (0,0) -- +(0,1);
\draw[near] (1,1) -- +(-1,0);
\draw[near] (1,1) -- +(0,-1);
\draw[near] (1,0) -- +(1,0);
\draw[near] (2,0) -- +(1,0);
\node[blackvertex] at (0,1) {};
\node[blackvertex] at (1,0) {};
\node[blackvertex] at (2,0) {};
\node[blackvertex] at (3,0) {};
\node at (0.65, 0.35) {$\sigma$};
\node at (3.6,-0.5) {(1c)};
\end{tikzpicture}
\vspace{15pt}

\begin{tikzpicture}[scale=0.85]
\unbalancedsquare{2};
\ladderplus{3}
\draw[near] (2,0) -- +(1,0);
\draw[near] (0,1) -- +(1,0);
\node[blackvertex] at (0,1) {};
\node[blackvertex] at (1,1) {};
\node[blackvertex] at (2,0) {};
\node[blackvertex] at (3,0) {};
\node at (1.65, 0.35) {$\sigma$};
\node at (3.6,-0.5) {(2a)};
\end{tikzpicture}
\hspace{13pt}
%
\begin{tikzpicture}[scale=0.85]
\ladderplus{3}
\draw[far] (1,0) -- +(1,0);
\draw[far] (1,0) -- +(0,1);
\draw[near] (2,1) -- +(-1,0);
\draw[near] (2,1) -- +(0,-1);
\draw[near] (2,0) -- +(1,0);
\draw[near] (0,1) -- +(1,0);
\node[blackvertex] at (0,1) {};
\node[blackvertex] at (1,1) {};
\node[blackvertex] at (2,0) {};
\node[blackvertex] at (3,0) {};
\node at (1.65, 0.35) {$\sigma$};
\node at (3.6,-0.5) {(2b)};
\end{tikzpicture}
\hspace{13pt}
%
\begin{tikzpicture}[scale=0.85]
\ladderplus{3}
\draw[near] (1,0) -- +(1,0);
\draw[near] (1,0) -- +(0,1);
\draw[far] (2,1) -- +(-1,0);
\draw[far] (2,1) -- +(0,-1);
\draw[near] (2,0) -- +(1,0);
\draw[near] (0,1) -- +(1,0);
\node[blackvertex] at (0,1) {};
\node[blackvertex] at (1,1) {};
\node[blackvertex] at (2,0) {};
\node[blackvertex] at (3,0) {};
\node at (1.65, 0.35) {$\sigma$};
\node at (3.6,-0.5) {(2c)};
\end{tikzpicture}
\end{center}
\caption{If $G=C_{4k+5} \cpg K_2$ contains any of the above configurations ((1a),(1b),(1c) in top row and (2a), (2b), (2c) in bottom row), then it admits a $(13,4)$-threshold-coloring with respect to $\mathcal L$. Black vertices go in $V_0$, $\sigma$ is an unbalanced square, configurations are shown up to symmetry. }
\label{fig:mod5:configs}
\end{figure}

Let us first argue $G$ is $(13,4)$-threshold-colorable, if it contains any of the configurations from Figure~\ref{fig:mod5:configs}. (This includes the case $n=5$.)
In all of the cases we may construct a square-wave path, which is delayed by two, and apply Lemma \ref{zigzag} in just the same way as in the previous cases.
The set $A$ contains (among others) all vertices depicted in black.
Observe that in case of (1a), (1b)  and (1c) configurations the square-wave path gets delayed by 2 immediately right of $\sigma$,
and in cases (2a), (2b), and (2c) two delays of $1$ appear at both ends of a square-wave path.
We can henceforth assume that none of the configurations from Figure~\ref{fig:mod5:configs} appears in $G$.

By the subclaim, $G$ admits an unbalanced square, and, without loss of generality, we may assume that $\sigma_0$ is unbalanced.
By (i), at least one of peripheral edges of $\sigma_1$ and $\sigma_{n-1}$ is a near edge, say $\overline{e}_1 \in N$.
The edge $\overline{e}_1$ can be used to delay a square-wave path starting at $u_1$ by $1$ (as it can continue through $u_2$). 
Another delay by one would finish the argument. For this we will assume that $n \geq 9$.

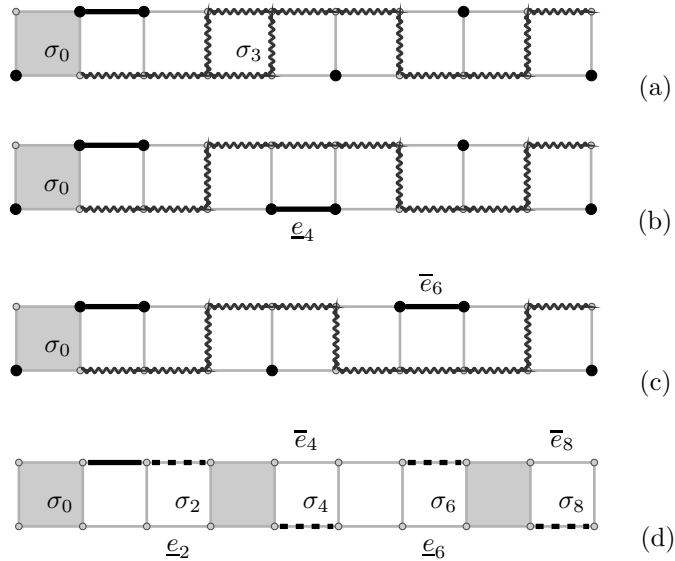
\begin{figure}
\begin{center}
\begin{tikzpicture}[scale=0.85]
\unbalancedsquare{1};
\ladder{9}
\draw[near] (1,1) -- +(1,0);
\draw[zigzag] (1,0) -- +(1,0);
\draw[zigzag] (2,0) -- +(1,0);
\draw[zigzag] (3,0) -- +(1,0);
\draw[zigzag] (3,1) -- +(1,0);
\draw[zigzag] (4,1) -- +(1,0);
\draw[zigzag] (5,1) -- +(1,0);
\draw[zigzag] (6,0) -- +(1,0);
\draw[zigzag] (7,0) -- +(1,0);
\draw[zigzag] (8,1) -- +(1,0);
\draw[zigzag] (3,0) -- +(0,1);
\draw[zigzag] (4,0) -- +(0,1);
\draw[zigzag] (6,0) -- +(0,1);
\draw[zigzag] (8,0) -- +(0,1);
\node at (0.65, 0.35) {$\sigma_0$};
\node at (3.65, 0.35) {$\sigma_3$};
\node[blackvertex] at (0,0) {};
\node[blackvertex] at (5,0) {};
\node[blackvertex] at (9,0) {};
\node[blackvertex] at (1,1) {};
\node[blackvertex] at (2,1) {};
\node[blackvertex] at (7,1) {};

\node at (10, -0.2) {(a)};

\end{tikzpicture}
\vspace{10pt}

\begin{tikzpicture}[scale=0.85]
\unbalancedsquare{1};
\ladder{9}
\draw[near] (1,1) -- +(1,0);
\draw[near] (4,0) -- +(1,0);
\draw[zigzag] (1,0) -- +(1,0);
\draw[zigzag] (2,0) -- +(1,0);
\draw[zigzag] (3,1) -- +(1,0);
\draw[zigzag] (4,1) -- +(1,0);
\draw[zigzag] (5,1) -- +(1,0);
\draw[zigzag] (6,0) -- +(1,0);
\draw[zigzag] (7,0) -- +(1,0);
\draw[zigzag] (8,1) -- +(1,0);
\draw[zigzag] (3,0) -- +(0,1);
\draw[zigzag] (6,0) -- +(0,1);
\draw[zigzag] (8,0) -- +(0,1);
\node at (0.65, 0.35) {$\sigma_0$};
\node at (4.5, -.35) {$\underline{e}_4$};
\node[blackvertex] at (0,0) {};
\node[blackvertex] at (4,0) {};
\node[blackvertex] at (5,0) {};
\node[blackvertex] at (9,0) {};
\node[blackvertex] at (1,1) {};
\node[blackvertex] at (2,1) {};
\node[blackvertex] at (7,1) {};
\node at (10, -0.2) {(b)};
\end{tikzpicture}
\vspace{5pt}

\begin{tikzpicture}[scale=0.85]
\unbalancedsquare{1};
\ladder{9}
\draw[near] (1,1) -- +(1,0);
\draw[near] (6,1) -- +(1,0);
\draw[zigzag] (1,0) -- +(1,0);
\draw[zigzag] (2,0) -- +(1,0);
\draw[zigzag] (3,1) -- +(1,0);
\draw[zigzag] (4,1) -- +(1,0);
\draw[zigzag] (5,0) -- +(1,0);
\draw[zigzag] (6,0) -- +(1,0);
\draw[zigzag] (7,0) -- +(1,0);
\draw[zigzag] (8,1) -- +(1,0);
\draw[zigzag] (3,0) -- +(0,1);
\draw[zigzag] (5,0) -- +(0,1);
\draw[zigzag] (8,0) -- +(0,1);
\node at (0.65, 0.35) {$\sigma_0$};
\node at (6.5, 1.35) {$\overline{e}_6$};
\node[blackvertex] at (0,0) {};
\node[blackvertex] at (4,0) {};
\node[blackvertex] at (9,0) {};
\node[blackvertex] at (1,1) {};
\node[blackvertex] at (2,1) {};
\node[blackvertex] at (6,1) {};
\node[blackvertex] at (7,1) {};

\node at (10, -0.2) {(c)};
\end{tikzpicture}
\vspace{5pt}

\begin{tikzpicture}[scale=0.85]
\unbalancedsquare{1};
\unbalancedsquare{4};
\unbalancedsquare{8};
\ladder{9}
\draw[near] (1,1) -- +(1,0);
\draw[far] (2,1) -- +(1,0);
\draw[far] (6,1) -- +(1,0);
\draw[far] (4,0) -- +(1,0);
\draw[far] (8,0) -- +(1,0);
\node at (0.65, 0.35) {$\sigma_0$};
\node at (2.65, 0.35) {$\sigma_2$};
\node at (4.65, 0.35) {$\sigma_4$};
\node at (6.65, 0.35) {$\sigma_6$};
\node at (8.65, 0.35) {$\sigma_8$};
\node at (2.5, -0.35) {$\underline{e}_2$};
\node at (6.5, -0.35) {$\underline{e}_6$};
\node at (4.5, 1.35) {$\overline{e}_4$};
\node at (8.5, 1.35) {$\overline{e}_8$};
\node at (10, -0.2) {(d)};
\end{tikzpicture}
\end{center}
\caption{Changing phase of a square-wave path in $C_9 \cpg K_2$ (a), (b), (c), towards a useful cut (d).}
\label{fig:mod5long}
\end{figure}

If $\sigma_3$ or $\sigma_7$ is an even square, then a further delay by $1$ can be achieved, see Figure~\ref{fig:mod5long}(a), the case $\sigma_3$ is an even square is shown.
Similarly, if one of the following edges $\overline{e}_2, \overline{e}_6, \underline{e}_4, \underline{e}_8$ is a near edge.
Figures~\ref{fig:mod5long}(b) and (c) contains the arguments where $\underline{e}_4 \in N$ or $\overline{e}_6 \in N$. 

Hence we may assume that both $\sigma_3$ and $\sigma_7$ are odd and hence unbalanced squares,
and that the edges $\overline{e}_2, \overline{e}_6, \underline{e}_4, \underline{e}_8$ all belong to $F$.
As configuration (2a) is not present at either $\sigma_3$ or $\sigma_7$, at least one of $\underline{e}_2$ or $\overline{e}_4$ is far,
and at least one of $\underline{e}_6$ and $\overline{e}_8$ is a far edge, see Figure~\ref{fig:mod5long}(d).
This implies that one of $\sigma_2$ or $\sigma_4$ contains a half-cut and so does one of $\sigma_6$ or $\sigma_8$.
This is by (i) not possible.

This proves Claim 3 in Case 4, unless $n=5$. So let us consider $G=C_5 \cpg K_2$. We shall split the analysis in half.

{\bf Subcase 4.1.} $n=5$ and there is an unbalanced square incident with exactly one near peripheral edge.

\begin{figure}
\begin{center}
\begin{tikzpicture}[scale=1.1]
\unbalancedsquare{3};
\ladder{5}
\draw[dte] (0,0) --node[black, below] {$\underline{e}_3$} +(1,0);
\draw[far] (1,0) --node[black, below] {$\underline{e}_4$} +(1,0);
\draw[dte] (2,0) --node[black, below] {$\underline{e}_0$} +(1,0);
\draw[far] (3,0) --node[black, below] {$\underline{e}_1$} +(1,0);
\draw[dte] (4,0) --node[black, below] {$\underline{e}_2$} +(1,0);
\draw[dte] (0,1) --node[black, above] {$\overline{e}_3$} +(1,0);
\draw[far] (1,1) --node[black, above] {$\overline{e}_4$} +(1,0);
\draw[dte] (2,1) --node[black, above] {$\overline{e}_0$} +(1,0);
\draw[near] (3,1) --node[black, above] {$\overline{e}_1$} +(1,0);
\draw[dte] (4,1) --node[black, above] {$\overline{e}_2$} +(1,0);
\draw[dte] (0,0) --node[black, left] {$e'_3$} +(0,1);
\draw[dte] (1,0) --node[black, left] {$e'_4$} +(0,1);
\draw[dte] (2,0) --node[black, left] {$e'_0$} +(0,1);
\draw[dte] (3,0) --node[black, right] {$e'_1$} +(0,1);
\draw[dte] (4,0) --node[black, right] {$e'_2$} +(0,1);
\draw[dte] (5,0) --node[black, right] {$e'_3$} +(0,1);
\node at (2.65, 0.35) {$\sigma_0$};
\node at (3.5, -1) {(initial setting)};
\end{tikzpicture}
\hspace{5mm}
\begin{tikzpicture}[scale=1.1]
\unbalancedsquare{3};
\ladder{5}
\draw[far] (0,0) --node[black, below] {$(5)$} +(1,0);
\draw[far] (1,0) --node[black, below] {$(1)$} +(1,0);
\draw[near] (2,0) --node[black, below] {$(8)$} +(1,0);
\draw[far] (3,0) --node[black, below] {$(1)$} +(1,0);
\draw[near] (4,0) --node[black, below] {$(3)$} +(1,0);
\draw[near] (0,1) --node[black, above] {$(6)$} +(1,0);
\draw[far] (1,1) --node[black, above] {$(1)$} +(1,0);
\draw[far] (2,1) --node[black, above] {$(8)$} +(1,0);
\draw[near] (3,1) --node[black, above] {$(1)$} +(1,0);
\draw[far] (4,1) --node[black, above] {$(2)$} +(1,0);
\draw[near] (0,0) --node[black, left] {$(7)$} +(0,1);
\draw[near] (1,0) --node[black, left] {$(9)$} +(0,1);
\draw[near] (2,0) --node[black, left] {$(8)$} +(0,1);
\draw[near] (3,0) --node[black, right] {$(8)$} +(0,1);
\draw[near] (4,0) --node[black, right] {$(4)$} +(0,1);
\draw[near] (5,0) --node[black, right] {$(7)$} +(0,1);
\node at (2.65, 0.35) {$\sigma_0$};
\node at (3.5, -1) {(final labels)};
\end{tikzpicture}
\end{center}
\caption{$G=C_5 \cpg K_2$: Unbalanced square $\sigma_0$ incident with exactly one near peripheral edge.}
\label{fig:mod5:one}
\end{figure}
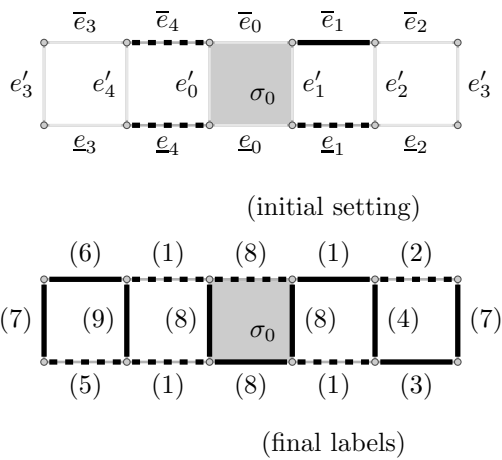

Let $\sigma_0$ be the unbalanced square incident with exactly one near peripheral edge $\overline{e}_1$.
The situation is depicted in Figure~\ref{fig:mod5:one}. The integer labels indicate the order in which we determine the labels of edges.
Here is how we infer:
\renewcommand\labelenumi{\theenumi}
\renewcommand{\theenumi}{(\arabic{enumi})}
\begin{enumerate}
\item $\overline{e}_1$ is a near edge and $\underline{e}_4, \underline{e}_1, \overline{e}_4$ are far edges, by assumption.
Note that $\sigma_4$ contains a half-cut, hence no other half-cuts are present by (i).
\item $\overline{e}_2$ is a far 
edge, as configuration (1a) from Figure~\ref{fig:mod5:configs} is not present,
\item $\sigma_2$ does not contain a half-cut by (i), hence $\underline{e}_2$ is a near edge,
\item similarly, $e'_2$ is a near edge, as otherwise a half-cut crosses squares $\sigma_1$ and $\sigma_2$.
\item Irrespective of the label of $e'_3$ we know that none of configurations (2a), (2b) is present around $\sigma_2$. This implies that $\underline{e}_3$ is a far edge.
\item Now $\overline{e}_3$ is near, otherwise also $\sigma_3$ contains a half-cut.
\item So is $e'_3$, as otherwise a half-cut lies across $\sigma_2 \cup \sigma_3$.
\item As configuration (1a) is not present at $\sigma_2$, the edge $\overline{e}_0$ is far, $e'_0$ is near by (ii) and not both of $\underline{e}_0, e'_1$ are far, again by (ii).
As $\sigma_0$ is unbalanced all of $e'_0,\underline{e}_0, e'_1$ are near edges.
\item Finally $e_4$ is a near edge by (ii).
\end{enumerate}  
Observe that in this case near edges span a path, and we can $(5,1)$-threshold-color $G$ with respect to $(N,F)$ by choosing colors $0,1,1,2,2,3,3,4,4,5$ along the ``near'' path starting at $v_0$.
As $(5,1) \leq (13,4)$, Claim 3 is settled in this case.

{\bf Subcase 4.2.} $n=5$ and no unbalanced square is incident with exactly one near peripheral edge.

We still have an unbalanced square $\sigma_0$, and it is incident with at least two near edges, one of them being $\overline{e}_1$.
Note that as configuration (2a) from Figure~\ref{fig:mod5:configs} is not present, the edge $\underline{e}_4$ is a far edge.
Now if $\underline{e}_1$ is a near edge, by a similar argument $\overline{e}_4$ is far, hence $\sigma_4$ contains a half-cut.
By the absence of configuration (1a) also $\underline{e}_2$ and $\overline{e}_2$ are far.
This implies that $\sigma_2$ contains another half-cut, which contradicts (i).

\begin{figure}
\begin{center}
\begin{tikzpicture}[scale=1.1]
\unbalancedsquare{3};
\ladder{5}
\draw[dte] (0,0) --node[black, below] {$\underline{e}_3$} +(1,0);
\draw[far] (1,0) --node[black, below] {$\underline{e}_4$} +(1,0);
\draw[dte] (2,0) --node[black, below] {$\underline{e}_0$} +(1,0);
\draw[far] (3,0) --node[black, below] {$\underline{e}_1$} +(1,0);
\draw[dte] (4,0) --node[black, below] {$\underline{e}_2$} +(1,0);
\draw[dte] (0,1) --node[black, above] {$\overline{e}_3$} +(1,0);
\draw[near] (1,1) --node[black, above] {$\overline{e}_4$} +(1,0);
\draw[dte] (2,1) --node[black, above] {$\overline{e}_0$} +(1,0);
\draw[near] (3,1) --node[black, above] {$\overline{e}_1$} +(1,0);
\draw[dte] (4,1) --node[black, above] {$\overline{e}_2$} +(1,0);
\draw[dte] (0,0) --node[black, left] {$e'_3$} +(0,1);
\draw[dte] (1,0) --node[black, left] {$e'_4$} +(0,1);
\draw[dte] (2,0) --node[black, left] {$e'_0$} +(0,1);
\draw[dte] (3,0) --node[black, right] {$e'_1$} +(0,1);
\draw[dte] (4,0) --node[black, right] {$e'_2$} +(0,1);
\draw[dte] (5,0) --node[black, right] {$e'_3$} +(0,1);
\node at (2.65, 0.35) {$\sigma_0$};
\node at (3.5, -1) {(initial setting)};
\end{tikzpicture}
\hspace{5mm}
\begin{tikzpicture}[scale=1.1]
\unbalancedsquare{3};
\unbalancedsquare{1};
\ladder{5}

\draw[far] (0,0) --node[black, below] {$(4)$} +(1,0);
\draw[far] (1,0) --node[black, below] {$(1)$} +(1,0);
\draw[far] (3,0) --node[black, below] {$(1)$} +(1,0);
\draw[near] (4,0) --node[black, below] {$(3)$} +(1,0);
\draw[far] (0,1) --node[black, above] {$(2)$} +(1,0);
\draw[near] (1,1) --node[black, above] {$(1)$} +(1,0);
\draw[near] (3,1) --node[black, above] {$(1)$} +(1,0);
\draw[far] (4,1) --node[black, above] {$(2)$} +(1,0);
\draw[far] (0,0) --node[black, left] {$(5)$} +(0,1);
\draw[near] (1,0) --node[black, right] {$(6)$} +(0,1);
\draw[near] (4,0) --node[black, right] {$(3)$} +(0,1);
\draw[far] (5,0) --node[black, right] {$(5)$} +(0,1);
\node at (2.65, 0.35) {$\sigma_0$};
\node at (0.65, 0.35) {$\sigma_3$};
\node at (3.5, -1) {(final labels)};
\end{tikzpicture}
\end{center}
\caption{$G=C_5 \cpg K_2$: Unbalanced square $\sigma_0$ incident with exactly two near peripheral edges.}
\label{fig:mod5:two}
\end{figure}
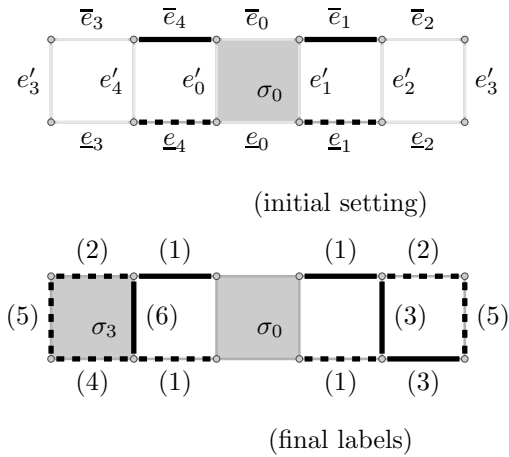

Therefore we may assume that $\overline{e}_4$ and $\overline{e}_1$ are the only near edges incident with $\sigma_0$. The situation is depicted in Figure~\ref{fig:mod5:two}

What follows is again a series of arguments, that will determine the next labels according to the integral labels in the same picture.
\begin{enumerate}
\item The labels of $\overline{e}_4, \underline{e}_4, \overline{e}_1, \underline{e}_1$ are determined by assumption. 
\item As configuration (1a) is absent, the edges $\overline{e}_2$ and $\overline{e}_3$ are both far edges.
\item Next consider edges $\underline{e}_3, \underline{e}_2, e'_4, e'_2$.
  At most one of them is a far edge, as otherwise we obtain either a pair of half-cuts or a $3$-cut.
  Neither is possible by (i) and (ii).
  We may by symmetry assume that both $\underline{e}_2, e'_2$ are near edges, and also that at least one of $\underline{e}_3, e'_4$ is a near edge. 
\item $\underline{e}_3$ is a far edge, otherwise one of the configurations (2a) or (2c) is present in square $\sigma_2$. 
\item $e'_3$ is a far edge, otherwise $\sigma_2$ is an unbalanced square incident with exactly one near peripheral edge, which is a situation settled in Subcase 4.1. 
\item By (ii), $e'_4$ is a near edge, which makes $\sigma_3$ an unbalanced square.
\end{enumerate}
But now a configuration (2a) from Figure~\ref{fig:mod5:configs} is present at $\sigma_3$ which is imposible.

This settles Claim 3 in the remaining case, and thus, in general.


{\bf Claim 4.}
Suppose that $G$ and $(N,F)$ do not admit a balanced parallel square with far spokes unless $n=3$.
Then $G$ is $(31,4)$-threshold-colorable with respect to $(N,F)$.

Let us construct a near-far-labeling $(N',F')$ such that $G$ and $(N',F')$ satisfy (ii), and (iii) of Claim 3.
Let $U$ be the set of vertices not incident with any edge from $N$ and take a maximal independent set $W \subseteq U$.
Let $D$ be the set of edges incident with some vertex of $W$. Then each vertex from $U$ is incident with an edge from $D \subseteq F$.
One readily checks that $G,(N',F'):=(N \cup D,F \setminus D)$ satisfy (ii) in Claim 3. Suppose $n \geq 4$ and assume, to the
contrary, that there is a balanced parallel square $\sigma$ with far spokes, all with respect to $(N',F')$; as there is no
such square with respect to the near-far-labeling $(N,F)$ by the conditions to $G,(N,F)$ of the Claim,
one of the peripheral edges of $\sigma$ is contained in $D$ (and has been moved from $F$ to $N'$);
but then one of its endvertices, say, $x$, is from $W$, so that {\em all} edges incident with $x$ have been moved from $F$ to $N'$, including one of the
spokes of $\sigma$, a contradiction. Therefore, $G,(N',F')$ satisfy (iii) of Claim 3.

If $G,(N',F')$ satisfy (i) of Claim 3, then, by Claim 3, there exists a $(13,4)$-threshold-coloring with respect to $G$;
otherwise, if (i) does not hold, then there exists an $(11,1)$-threshold-coloring of $G$ with respect to $(N',F')$ by Claim 2.
As $(13,4) \leq (26,4)$ and $(11,1) \leq (26,4)$, we find, in either case,
a $(26,4)$-threshold-coloring $c'$ of $G$ with respect to $(N',F')$. Setting $c(x):=c'(x)$ for $x \in V(G) \setminus W$ and $c(x)=30$ for $x \in W$
produces an $(31,4)$-threshold-coloring with respect to $(N,F)$, proving Claim 4.


We now finish by proving by induction on $n$ that for every graph $G=C_n \cpg K_2$, $n \geq 3$, and every near-far-labeling $(N,F)$ of $G$
there exists an $(31,4)$-threshold-coloring.
The induction starts for $n=3$ by Claim 4. Let us assume $n \geq 4$ and consider $G=C_n \cpg K_2$.
If $G$ admits no balanced parallel square with far spokes then again the statement follows by Claim 4.
Hence we may assume that $G$ admits a balanced square $\sigma$ with far spokes.
In this case we may contract the (near) peripheral edges of such a square and lift the coloring of the result back to $G$. In details:
Let $\{u_1 u_2, v_1 v_2 \} = E(\sigma) \cap N$.
Let $G'$ defined to be the prism obtained by contracting both near edges in $\sigma$,
whose edge labeling $\pi' = (N \setminus \{u_1 u_2, v_1 v_2 \}, F)$ matches the one of $G$.
By induction, there exists an $(31,4)$-threshold-coloring of $G'$ with respect to $\pi'$.
If $v,u$ are vertices obtained by contracting edges $v_1 v_2$ and $u_1 u_2$, respectively,
then by setting $c(u_1)=c(u_2) = c(u)$ and $c(v_1)=c(v_2) = c(v)$ and $c(x):=c'(x)$ for all $x \in V(G) \setminus \{u_1,u_2,v_1,v_2\}$
we obtain an $(31,4)$-threshold-coloring of $G$, as promised.
\hspace*{\fill}$\Box$

\section{The Petersen graph}

This entire section will show how to threshold color the {\sc Petersen} graph, which is denoted by $G$ throughout.
In order to keep the arguments as clear as possible, we will not try to make an effort in estimating the coloring parameters. 

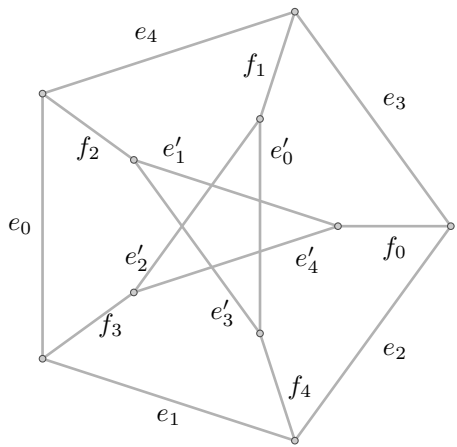
\begin{figure}
\begin{center}
\begin{tikzpicture}[scale=1.5]
\petersencoords;
\petersenlabelededges;
\petersennodes;
\end{tikzpicture}
\end{center}
\caption{The {\sc Petersen} graph $G$ and its edge labels.}
\label{fig:petersen:labels}
\end{figure}

The automorphism group of $G$ is isomorphic to $S_5$.
Apart from being vertex- and edge-transitive, $G$ is also $3$-arc transitive, that is,
every directed path of three edges can be mapped to every other such path using an automorphism of $G$.
Our analysis of threshold colorings of $G$ will use notation from Figure~\ref{fig:petersen:labels}.  

 
We shall first consider a pair of structures that will enable us to threshold color $G$ with respect to a particular labeling. 

Given an edge labeling $\pi:=(N,F)$, an \emph{$xxyyzz$-cycle} is a cycle of length $6$ in $G$ along edges
$\varepsilon_1$, $\varepsilon_2$, $\varepsilon_3$, $\varepsilon_4$, $\varepsilon_5$, $\varepsilon_6$,
so that for every $i \in \{1,2,3\}$ the edges $\varepsilon_{2i-1}$ and $\varepsilon_{2i}$ both belong to $N$ or both belong to $F$.
\begin{lem}
\label{lem:pet:6cycle}
Assume that $G,\pi$ admit an $xxyyzz$-cycle $C$. Then $G$ is $(5,1)$-threshold-colorable with respect to $\pi$.
\end{lem}

{\bf Proof.}
Assuming the above notation, let $u_1, u_2, u_3$ be the common vertices of $\varepsilon_6, \varepsilon_1$; $\varepsilon_2,\varepsilon_3$; and $\varepsilon_4,\varepsilon_5$, respectively. 
Set $A = \{u_1,u_2,u_3\}$, and $B = V(G) \setminus A$, and $M=\emptyset$.
We shall argue that we can use Lemma~\ref{zigzag} in this case.

Note first that $A$ is an independent set, hence (i) of Lemma \ref{zigzag} holds, and that $G[B]$ is a tree (indeed a claw $K_{1,3}$), which implies (ii) of Lemma \ref{zigzag}.
Condition (iii) of Lemma \ref{zigzag} is trivially true.
By the label structure of $C$ we also have (iv) of Lemma \ref{zigzag}, as every vertex $v \in B$ which is adjacent to (exactly) two vertices of $A$ uses a pair of edges with the same label.
\hspace*{\fill}$\Box$

A similar but a bit trickier argument deals with a specific $5$-cycle. Again, let us fix a labeling $\pi=(N,F)$.
A $5$-cycle $C$ along edges $\varepsilon_0,\varepsilon_1,\varepsilon_2,\varepsilon_3,\varepsilon_4$
is an {\em $Nxxyy$-cycle} if $\varepsilon_0 \in N$, and for every $i \in \{1,2\}$ the edges $\varepsilon_{2i-1}$ and $\varepsilon_{2i}$ are either both in $N$ or both in $F$.

\begin{lem}
\label{lem:pet:5cycle}
Assume that $G,\pi$ admit an $Nxxyy$-cycle $C$. Then $G$ is $(14,4)$-threshold-colorable with respect to $\pi$.
\end{lem}

{\bf Proof.}
Observe that $(5,1) \le (14,4)$, so constructing a $(5,1)$-threshold-coloring suffices.

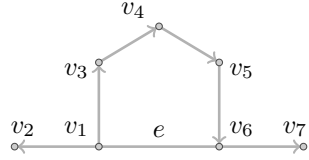
\begin{figure}
\begin{center}
\begin{tikzpicture}[scale=0.8]
\redpetersen;
\end{tikzpicture}
\end{center}
\caption{$G[V_1]$ from the proof of Lemma~\ref{lem:pet:5cycle}.}
\label{fig:redpetersen}
\end{figure}

Let $A$ be the vertex set containing both end vertices of $\varepsilon_0$ and the common end vertex of $\varepsilon_2$ and $\varepsilon_3$, and let $B = V(G) \setminus A$.
Note that $B$ induces the subgraph shown in Figure~\ref{fig:redpetersen}; we will use the notation from the very same figure.
 
Let us construct a $(5,1)$-threshold-coloring $c$ of $G-e$.
Set $c(x) = 0$ for every $x \in A$ and then let us color the remaining vertices $v_1, \ldots, v_7$ in order of their indices.
Setting $c(v_1)=1$ is admissible, as $v_1$ has no neighbor in $V_0$.
Next let us inductively choose colors of $c(v_2), c(v_3), c(v_4), c(v_5) \in \{-2,-1,1,2\}$ as in the proof of Lemma~\ref{zigzag}. 

Next also $v_6$ has no neighbor in $V_0$, so we may choose $c(v_6) \in \{-1,2\}$ so that the difference $|c(v_5) - c(v_6)|$ matches the label of $v_5v_6$
(again, the sign changes if and only if $v_5v_6 \in F$), and finally extend the coloring to $c(v_7)$ as usual: the --- by construction constant --- label of the $v_7,A$-edges
determines the absolute value, and the label of $v_6v_7$ determines the sign of $c(v_7)$.

If the number of far edges in the only cycle of $G[B]$ is even, then the $|c(v_1)-c(v_6)|$ matches the label of $v_1v_6$,
as in the proof of Lemma \ref{zigzag}, and we have a $(5,1)$-threshold-coloring. 

In case that the number of far edges along the cycle of $G[B]$ is odd, we may encounter two possibilities.
Recall that by construction $c(v_1)=1$.
Either $v_1v_6 \in N$ and  $c(v_6)=-1$ or $v_1v_6 \in F$ and $c(v_6)=2$.
In both cases we shall construct a $(14,4)$-threshold-coloring $c'$. For $x \in V(G) \setminus \{v_1,v_6\}$ let  $c'(x)=3 c(x)$.
In the first case let us set $c'(v_1)=2 = 3 c(v_1)-1$ and $c'(v_6)=-2 = 3 c(v_6)+1$ and in the second case we set $c'(v_1)=2 = 3 c(v_1)-1$ and $c'(v_6)=7= 3 c(v_6)+1$.
This implies that, in either case, $c'$ is a $(14,4)$-threshold-coloring of $G$ with respect to $\pi$.
\hspace*{\fill}$\Box$

Note that Lemma~\ref{zigzag} succeeds in constructing a $(13,4)$-threshold-coloring, and Lemma~\ref{lem:pet:5cycle} requires an additional color in some cases
($7$ has been assigned in the final paragraph of the preceeding proof).


Let us call a cut $S \subseteq F$ a {\em far-edge-cut}.
In what follows we shall prove that every edge labeling $\pi = (N,F)$ of $G$
admits either a far-edge-cut, a $Nxxyy$-cycle or a $xxyyzz$-cycle.

We may take care of far-edge-cuts inductively by assuming that $G$ is threshold colorable with respect to every edge labeling $\pi'= (N',F')$ for which $|F'| < |F|$.
The induction basis $F' = \emptyset$ is trivial, as a constant function can serve as a threshold coloring with threshold $0$ in this case.
If $S \subseteq F$ is a (minimal) far-edge-cut we can easily construct a threshold coloring of $G$ with respect to $\pi$ from a threshold coloring of $G$
with respect to $(N \cup S,F \setminus S)$ by increasing the colors simultaneously in one component of $G-S$ by some moderate constant.

In what follows we shall assume that $G,\pi$ contains no $xxyyzz$-cycles, no $Nxxyy$-cycles, and also no far-edge-cuts.
The latter implies that $(V(G),N)$ is connected; in particular,
(i) every vertex is incident with at least one near edge, and 
(ii) $|N| \geq 9$, so that $|F| \leq 6$.

Given a near-far-labeling $\pi$ we shall split our analysis according to the length of a longest ``far-edge-path''.

{\bf Case 1.} $G$ contains a path $P$ with $E(P) \subseteq F$ of length at least $4$. 

There are two subcases to consider.
The end vertices of the initial length $4$ segment of $P$ may be adjacent ($P$ is a subpath of some $5$-cycle), or else $P$ is a segment of a $6$-cycle.
Let us consider the former option first. 

\begin{figure}
\begin{center}
\begin{tikzpicture}[scale=0.6]
\petersencoords;
\petersenedges;
\foreach \x/\y in {p0/p1, p2/p3, p3/p4, p0/p4, p5/p8}{
  \draw[far] (\x) -- (\y);
}
\foreach \x/\y in {p3/p8, p4/p9, p0/p5}{
  \draw[near] (\x) -- (\y);
}
\petersennodes;
\end{tikzpicture}
\hspace{2cm}
\begin{tikzpicture}[scale=0.6]
\petersencoords;
\petersenedges;
\foreach \x/\y in {p2/p3, p3/p4, p0/p4, p0/p5}{
  \draw[far] (\x) -- (\y);
}
\foreach \x/\y in {p0/p1, p2/p1, p3/p8, p4/p9, p5/p8, p6/p8, p5/p7}{
  \draw[near] (\x) -- (\y);
}
\petersennodes;
\end{tikzpicture}
\end{center}
\caption{Four far edges on a path: either lying on a $5$-cycle (left) or not lying on a $5$-cycle (right) case.}
\label{fig:4:c1}
\end{figure}
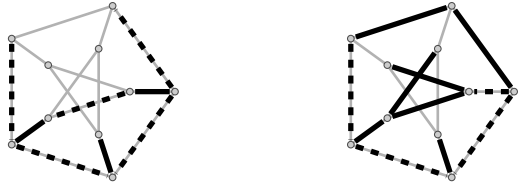

We may without loss of generality assume that $e_0, e_1, e_2, e_3$ are far edges.
This implies that $e_4$ is a far edge, as otherwise we have a $Nxxyy$-cycle.
By the absence of far-edge-cuts we have $f_3, f_4 ,f_0 \in N$.
Now $e'_4 \in F$ as otherwise we have a $Nxxyy$-cycle along edges $e'_4, f_3, e_1, e_2, f_0$, see Figure~\ref{fig:4:c1}(left).
By a similar argument at least one of $f_1, e'_0$ is a far edge, and also at least one of $f_2, e'_3$ is a far edge.
This implies that at least $7$ edges belong to $N$ which is absurd. 

Let us next assume that the four far edges of $P$ are $e_0, e_1, e_2, f_0$ and,
excluding the previous case, that no $5$-cycle contains four far edges.
Hence $e_3, e_4 \in N$, and by the absence of $3$-cuts we have $f_3,f_4 \in N$.
By excluding the previous case once more the edge $e'_4$ is a near edge.
As there are no $xxyyzz$-cycles, at least one of $f_1, e'_0$ is a far edge,
and also at least one of $f_2, e'_3$ is a far edge as there are no $Nxxyy$-cycles. 
Since at most $6$ edges are far, both $e'_2$ and $e'_1$ are near edges, see Figure~\ref{fig:4:c1}(right).
Now if any of the edges $f_1,f_2$ is a near edge, we obtain a $xxyyzz$-cycle.
Hence both are far, and consequently the edges $e'_3$ and $e'_0$ from the interior $5$-cycle are both near.
This is not possible as a $5$-cycle whose edges are uniformly near (here: the ``inner'' cycle) is a $Nxxyy$-cycle.   

{\bf Case 2.} A longest path $P$ in $G$ with $E(P) \subseteq F$ contains exactly $3$ edges.

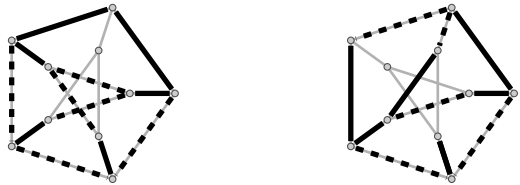
\begin{figure}
\begin{center}
\begin{tikzpicture}[scale=0.6]
\petersencoords;
\petersenedges;
\foreach \x/\y in {p2/p3, p3/p4, p0/p4, p5/p8, p7/p9, p5/p7}{
  \draw[far] (\x) -- (\y);
}
\foreach \x/\y in {p0/p1, p1/p2, p2/p7, p0/p5, p3/p8, p4/p9}{
  \draw[near] (\x) -- (\y);
}
\petersennodes;
\end{tikzpicture}
\hspace{2cm}
\begin{tikzpicture}[scale=0.6]
\petersencoords;
\petersenedges;
\foreach \x/\y in {p3/p4, p0/p4, p1/p2, p5/p8, p1/p6}{
  \draw[far] (\x) -- (\y);
}
\foreach \x/\y in {p2/p3, p0/p1, p3/p8, p4/p9, p4/p9, p0/p5, p6/p8}{
  \draw[near] (\x) -- (\y);
}
\petersennodes;
\end{tikzpicture}
\end{center}
\caption{Three far edges on a path (left) and at most two far edges on a path case (right).}
\label{fig:3:c1}
\end{figure}

Without loss of generality let $e_0, e_1, e_2$ be far edges.
By excluding previous cases and far-edge-cuts the edges $e_3, e_4, f_2, f_0, f_3, f_4$ are all near edges.
As we have no $Nxxyy$-cycles the edges $e'_4$ and $e'_3$ are both far, and so is $e'_1$, see Figure~\ref{fig:3:c1}(left).
By the far-edge-count alone we infer that $f_1, e'_0$ (and $e_2'$) are near,
and we have an $xxyyzz$-cycle along edges $e_4, f_1, e'_0, f_4, e_1, e_0$, which is absurd.

{\bf Case 3.} A longest path $P$ in $G$ with $E(P) \subseteq F$ contains exactly $2$ edges.

Again, we may without loss of generality assume that $e_1, e_2$ are far edges,
which in turn implies that $e_0, e_3, f_0, f_3$ are near (by maximality of an $F$-path) and $f_4$ is near (by the absence of far-edge-cuts).
By the absence of $Nxxyy$-cycles we infer that both $e_4$ and $e'_4$ are far edges. 

Next let us consider a pair of edges $e'_2, f_1$.
As the cycle along $e_3, e_2, e_1, f_3, e'_2, f_1$ is not a $xxyyzz$-cycle, at least one of $e'_2, f_1$ is a far edge,
and since there is no path of $3$ far edges, exactly one of them is.
Focusing on the cycle along $e_3, f_0, e'_4, e'_2, f_1$, which cannot be an $Nxxyy$-cycle, implies that $e'_2 \in N$ and $f_1 \in F$,
but in this case the cycle along $e_0, f_3, e'_2, f_1, e_4$ is a $Nxxyy$-cycle, see Figure~\ref{fig:3:c1}(right).

We are left with the final case.

{\bf Case 4.} $F$ forms a matching.

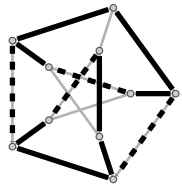
\begin{figure}
\begin{center}
\begin{tikzpicture}[scale=0.6]
\petersencoords;
\petersenedges;
\foreach \x/\y in {p0/p4, p2/p3, p5/p7, p6/p8}{ 
  \draw[far] (\x) -- (\y);
}
\foreach \x/\y in {p3/p4, p4/p9, p0/p1, p0/p5, p1/p2, p2/p7, p3/p8, p6/p9}{
  \draw[near] (\x) -- (\y);
}
\petersennodes;
\end{tikzpicture}
\end{center}
\caption{Far edges form a matching.}
\label{fig:c:match}
\end{figure}

Every matching of size $5$ in the {\sc Petersen} graph is an edge cut, 
so we may assume that $|F| \leq 4$. 
Again we may assume that $e_2$ is a far edge, which in turn makes $e_1, e_3, f_4, f_0 \in N$. 
By the absence of a $xxyyzz$-cycle, we may without loss of generality assume that $e_0 \in F$ and $e_4 \in N$, and consequently also $f_2, f_3 \in N$.
Now by the absence of an $Nxxyy$-cycle, $e'_1$ is a far edge, and also exactly one of $e'_2, e'_0$ (by symmetry we may assume the former), see Figure~\ref{fig:c:match}.
This makes the cycle along $e_4, f_1, e'_0, e'_3, f_2$ an $Nxxyy$-cycle. 


{\bf Addresses of the authors.}

\parbox{4cm}{
{\sc Ga\v{s}per Fijav\v{z}} \\
Faculty of Computer and Information Science\\
University of Ljubljana \\
Ve\v{c}na pot 113\\
1000 Ljubljana \\
Slovenia}
\hspace*{\fill}
\parbox{6cm}{
{\sc Matthias Kriesell} \\
Institut f\"ur Mathematik\\
Technische Universit\"at Ilmenau \\
Weimarer Stra{\ss}e 25 \\
98693 Ilmenau \\
Germany}

\end{document}